\documentclass[moor]{informs3} 

\usepackage{natbib}
 \NatBibNumeric
 \def\bibfont{\small}%
 \bibpunct[, ]{[}{]}{,}{n}{}{,}%

\usepackage[colorlinks=true,breaklinks=true,bookmarks=true,urlcolor=blue, 
     citecolor=blue,linkcolor=blue,bookmarksopen=false,draft=false]{hyperref}

\def\EMAIL#1{\href{mailto:#1}{#1}}

\TheoremsNumberedThrough
\EquationsNumberedThrough

\usepackage{theorem}
\usepackage{epic,eepic}
\usepackage{graphicx}
\usepackage{algorithm}
\usepackage{amssymb,amsmath,amsfonts,mathrsfs}
\usepackage{enumitem}
\usepackage{etoolbox,verbatim,color}
\usepackage{microtype}
\usepackage{url}
\usepackage{algorithmic}
\usepackage{cleveref}
\newcommand{\eqa}{\stackrel{\rm(a)}{=}}
\newcommand{\eqb}{\stackrel{\rm(b)}{=}}
\newcommand{\eqc}{\stackrel{\rm(c)}{\equiv}}

\newcommand{\mADMM}{\hyperref[alg:mADMM]{{\sf mADMM}}}

\newcommand{\mL}{\mathcal{L}}
\newcommand{\mB}{\mathcal{B}}

\newcommand{\bbR}{\mathbb R}

\newcommand{\QED}{\hfill $\square$}

\newcommand{\prox}{\ensuremath{\operatorname{Prox}}}

\DeclareMathOperator*{\minimize}{minimize}

\newtheorem{condition}{\textbf{Condition}}
\newcommand{\mcalB}{\mathcal{B}}
\begin{document}

 \TITLE{Multiblock ADMM for nonsmooth nonconvex optimization with nonlinear coupling constraints}

\ARTICLEAUTHORS{%
\AUTHOR{Le Thi Khanh Hien}
\AFF{Department of Mathematics and Operational Research, University of Mons, Belgium, \EMAIL{thikhanhhien.le@umons.ac.be}}
\AUTHOR{Dimitri Papadimitriou}
\AFF{3NLab, Huawei Belgium Research Center (BeRC)\\ \EMAIL{ dimitrios.papadimitriou.ext@huawei.com} }
}

	\ABSTRACT{		This paper proposes a multiblock alternating direction method of multipliers for solving a class of multiblock nonsmooth nonconvex optimization problem with nonlinear coupling constraints. We employ a majorization minimization procedure in the update of each block of the primal variables. Subsequential and global convergence of the generated sequence to a critical point of the augmented Lagrangian are proved.  We also establish iteration complexity and provide preliminary numerical results for the proposed algorithm.
	}		
  		
  	
\maketitle

\section{Introduction}
\label{sec:intro}
We consider the following multiblock  optimization problem with nonlinear coupling constraints 
\vspace{-0.1in} 
\begin{equation}
\label{eq:model}
  \minimize_{x,y}  \; F(x_1,\ldots,x_m) + h(y) 
 \quad 
\text{ s.t}  \quad \phi(x) + \mathcal B y = 0,
\end{equation} 
where $x$ can be decomposed into $m$ blocks $ x=(x_1,\ldots,x_m)$ with $x_i\in \mathbb R^{ n_i}$, $ n=\sum_{i=1}^m  n_i$, $y \in \mathbb R^{ q}$, $\phi$ is a nonlinear mapping from $\mathbb R^{n}$ to $\mathbb R^{ s}$ defined by $\phi(x)=(\phi_1(x),\ldots,\phi_s(x))$, $\mathcal B$ is a linear map,  $F(x)=f(x)+ g(x)$,   $g(x)=\sum_{i=1}^m g_i(x_i)$, where $f: \mathbb R^{n} \to \mathbb R$ is a continuous (but possibly nonsmooth) function,  $g_i:\mathbb R^{ n_i} \to \mathbb R\cup \{+\infty\}$ are  proper lower semi-continuous (lsc) functions for $i=1,2,\dots,m$, and $h:\mathbb R^{ q} \to \mathbb R$ is a  differentiable function. Note that $F$ and $h$ can be nonconvex. Throughout the paper, we assume the following assumptions.
\begin{assumption}
(A1)  $F$ satisfies 
$\partial F(x) = \partial_{x_1} F(x) \times \ldots \times \partial_{x_m} F(x)$, where $\partial F$ denote the limiting  subdifferential of $F$ (see  \cite[Definition 8.3]{RockWets98} for its definition).

(A2)  $\nabla h$ is $L_h$-Lipschitz continuous.

(A3) $\sigma_\mB:=\lambda_{\min} (\mB \mB^*)>0$ ($\sigma_\mcalB$ is the smallest eigenvalue of $\mB\mB^*$).

(A4) $F(x) + h(y)$ is lower bounded. 
\end{assumption}
(Note that Assumption (A1) is satisfied when $f$ is a sum of a continuously differentiable function and a block separable function;  see~\cite[Proposition 2.1]{Attouch2010}.) 
While the Alternating Direction Method of Multipliers (ADMM) has gained significant attention and proven highly effective in addressing multiblock composite optimization problems with linear constraints, see e.g., \cite{Hien2021_iADMM,Wang2019}, the exploration of ADMM's applicability to nonconvex nonsmooth problems with  \emph{nonlinear coupling constraints} remains relatively limited. 
    The authors in \cite{Bolte2018_ADMM} considers Problem~\eqref{eq:model} with $m=1$ and a general nonlinear coupling constraint. They  propose a universal framework to study  global convergence analysis of Lagrangian sequences, introduce the notion of information zone, and propose an adaptive regime to detect this zone in finitely many steps and force the iterates to remain in the zone. Note that the essential assumption for ensuring the information zone's role and the global convergence of the Lagrangian sequence in the framework of \cite{Bolte2018_ADMM} is the boundedness of the multiplier sequence, as depicted in \cite[Lemma 1]{Bolte2018_ADMM}. Moreover, the choice of parameters within the context of \cite{Bolte2018_ADMM} is closely tied to the upper bound constant of the multiplier sequence, as discussed in \cite[Remark 7]{Bolte2018_ADMM}. In practice, determining this upper bound for the multiplier can be an exceedingly challenging, if not impossible, task.  
    The authors in \cite{Cohen2021} consider Problem~\eqref{eq:model} with $m=1$. They design a proximal linearized ADMM in which the proximal parameter is generated dynamically by a backtracking procedure. Although the boundedness assumption of the multiplier sequence is not required in \cite{Cohen2021}, the backtracking procedure used in \cite{Cohen2021}, on one hand, requires to repeatedly evaluate the values of $f$ and $\phi$ to guarantee a descent property in updating $x$, and on the other hand, relies on the boundedness assumption of the generated sequence 
to guarantee the boundedness of the proximal parameters produced by the backtracking procedure. However, \cite{Cohen2021} does not provide a sufficient condition to guarantee the important boundedness assumption. 
    
 In this paper, we propose a multiblock alternating direction method of multipliers (\mADMM) for solving Problem \eqref{eq:model}. Different from the general framework of \cite{Bolte2018_ADMM}, we specify the update of the primal \emph{block} variables by embedding majorization minimization (MM) procedure. MM step can recover proximal point, proximal gradient and mirror descent step by using suitable surrogate functions depending on the structure of the objective function, see \cite{Titan2020,Razaviyayn2013}. Hence, employing MM and conducting the convergence analysis based on the surrogate functions allow us to not only better explore the structure of the problem but also unify the convergence analysis of special cases.  We prove the subsequential convergence standard assumptions, and prove the global convergence with some additional assumptions. We also establish the iteration complexity and provide numerical results for \mADMM~(we note that iteration complexity was not established and numerical results were not presented in \cite{Bolte2018_ADMM} and \cite{Cohen2021}).

 The paper is organized as follows. In the next section, we provide some preliminary knowledge to support the forthcoming analysis.  In Section~\ref{sec:mADMM}, we describe \mADMM~and conduct the convergence analysis for it.  We also provide a sufficient condition such that the generated sequence by \mADMM\, is bounded. We report numerical results in Section~\ref{sec:numerical} and conclude the paper in Section~\ref{sec:conclusion}. 

\textbf{Notation.} We denote $[m]=\{1,\ldots,m\}$. For a mapping $\phi : \mathbb R^{n} \to \mathbb R^{s}$, we denote $\nabla \phi (x) =[\nabla \phi_1(x) \ldots \nabla \phi_s(x)]$ ($\nabla \phi (x)^\top$ would be the Jacobian), 
 and $\nabla_{x_i} \phi (x) =[\nabla_{x_i} \phi_1(x) \ldots \nabla_{x_i} \phi_s(x)]
 $.

 \section{Preliminaries}
\label{sec:problem}

\subsection{Augmented Lagrangian and $\varepsilon$-stationary point.} 
The augmented Lagrangian function of  Problem~\eqref{eq:model} is 
\vspace{-0.1in} 
\begin{equation}
\label{Lagrangian}
\begin{array}{ll}
&\mathcal L_\beta(x, y, \omega)=\varphi_\beta(x,y,\omega) + \sum_{i=1}^m g_i(x_i),\\
\text{where} &\varphi_\beta(x,y,\omega)= f(x) + h(y)
 +  \langle \omega,\phi(x) + \mcalB y \rangle + \frac{\beta}{2} \| \phi(x) +  \mcalB y\|^2.
\end{array} 
\end{equation} 
and $\beta>0$ is a penalty parameter.  
 Let $(x,y,\omega)$ be a critical point of $\mL_\beta$, that is 
\vspace{-0.1in} \begin{equation}
\label{cond:critical_point}
   \begin{array}{ll}
 &   0 \in \partial_{x_i} F + \nabla_{x_i} \phi(x) \big(\omega + \beta(  \phi( x) +\mcalB y )\big),\\
    & \nabla h(y) + \nabla \mcalB y \big(\omega + \beta(  \phi( x) +\mcalB y )\big)=0,  \quad 
    \phi( x) +\mcalB y=0.     
\end{array} 
\end{equation} 
The conditions in \eqref{cond:critical_point} imply that   
$    0 \in \partial_{x_i} F + \nabla_{x_i} \phi(x) \omega$, $ \nabla h(y) +   \mcalB^* \omega =0$,  and $
     \phi( x) +\mcalB y =0$.     
Hence $(x,y)$ is a stationary point of  Problem \eqref{eq:model}. Conversely, if $(x,y)$ is a stationary point of Problem \eqref{eq:model} then there exists $\omega$ such that $(x,y,\omega)$ satisfies~\eqref{cond:critical_point}, which means that $(x,y,\omega)$ is a critical point of $\mL_\beta$. Therefore, finding a stationary point of \eqref{eq:model} is equivalent to finding a critical point of the augmented Lagrangian $\mL_\beta$. 
\begin{definition}
\label{def:epsilon_opt}
We call $(x,y)$ an $\varepsilon$-stationary point of \eqref{eq:model} if there exists $\omega$ and $\chi_i\in \partial_{x_i} F(x)$ such that
$ R_i=\| \chi_i + \nabla_{x_i} \phi(x) \omega\| \leq \varepsilon$, $R_y=\| \nabla h(y) +  \mcalB^*\omega\|\leq \varepsilon$,   and $R_c =\| \phi( x) +\mcalB y\|\leq \varepsilon.
$
\end{definition}

\subsection{Block surrogate function.}
The notations used in this subsection are independent of the notations of other sections. 
Let us first remind the proximal gradient method (PG), see e.g.,  \cite{Boyd2014}, for solving the following composite optimization problem 
\begin{equation}
\label{eq:composite_prob}
\min_{x\in\mathbb R^n} f(x) + g(x),
\end{equation} 
where  $f$ is a continuous differentiable function, $\nabla f(x)$ is assumed to be $L$-Lipschitz continuous and $g$ is a lower semicontinuous function. The PG step is   
\begin{equation}
\label{proximal_point}
x^{k+1}\in\prox_{\frac{1}{L}g}(x^k-\frac{1}{L} \nabla f(x^k)),
\end{equation} 
where the proximal mapping is defined as  
$$\prox_{g} (x):=\arg\min_{y} g(y) + \frac{1}{2}\|y-x\|^2. 
$$ 
As $\nabla f$ is $L$-Lipschitz continuous, the descent lemma (see \cite{Nesterov2004}) gives us
\begin{equation}
\label{descentlemma}
 f(x) \leq f(x^k) + \langle \nabla f(x^k),x-x^k) + \frac{L}{2} \|x-x^k\|^2.
\end{equation} 
Denote the right hand side of \eqref{descentlemma} by $u(x,x^k)$. We notice that 
\begin{equation}
\label{surrogate_function}
u(x,x)=f(x) \quad \mbox{and} \quad  u(x,y)\geq f(x),\, \forall\, x,y.
\end{equation}
Each function $u(\cdot,\cdot)$ that satisfies \eqref{surrogate_function} is called a surrogate function of $f$. The PG step in \eqref{proximal_point} can be rewritten in an equivalent form 
\begin{equation}
\label{MM}
x^{k+1}=\arg\min_x u(x,x^k) + g(x),
\end{equation}
where $u$ is the Lipschitz gradient surrogate defined by the right hand side of \eqref{descentlemma}. By choosing suitable surrogate functions, many first order methods such as proximal methods, mirror descent methods, etc., can also be rewritten in the majorization minimization (MM) form. Convergence analysis of MM algorithm (which can be roughly described in \eqref{MM}) would unify the convergence analysis of the algorithms that correspond to different choices of the surrogates. 

When $x$ has multiple blocks $x=(x_1,\ldots,x_m)$, $x_i\in \mathbb R^{n_i}$, and we assume that $g(x)=\sum g_i(x_i)$, block coordinate descent (BCD) method  is a well-known approach to solve \eqref{eq:composite_prob}. BCD updates one block at a time while fixing the values of the other blocks. The condition in~\eqref{surrogate_function} is then extended to deal with functions that have multiblock variables as follows.
We adopt the definition of block surrogate functions from~\cite{Razaviyayn2013}. 
 
\begin{definition}[Block surrogate function]
\label{def:surrogate-block} Let $\mathcal X_i\subseteq \mathbb R^{n_i}$, $\mathcal X =\mathcal X_1 \times \ldots \times \mathcal X_m \subseteq \mathbb R^n$. 
A continuous function $u_i:\mathcal X_i \times \mathcal X \to \mathbb R  $ is called a block surrogate function of $f$ on $\mathcal X$ with respect to block $x_i$
if 
\vspace{-0.1in}$$ u_i(z_i,z) = f(z)\quad\mbox{and} \quad  
    u_i(x_i,z) \geq f(x_i,z_{\ne i})\,\,\mbox{for all} \,\, x_i\in\mathcal X_i, z\in \mathcal X, 
    $$
where  
    $(x_i,z_{\ne i})$ denotes $(z_1,\ldots,z_{i-1},x_i,z_{i+1},\ldots,z_m).$ 
The block approximation error is defined as 
$ e_i(x_i,z):=u_i(x_i,z) - f(x_i,z_{\ne i}).$
\end{definition}
When writing the MM step for the update of block $i$, cf. \eqref{eq:x_iupdate-2}, we can perceive that the current point is $z$ and the value $z_i$ of block $i$ will be updated by $\arg\min_{x_i} u_i(x_i,z_{\ne i}) + g_i (x_i)$.
\begin{example}
\label{example:surrogates}
\begin{enumerate}
    \item[(i)] A proximal surrogate (see e.g., \cite{Attouch2009,Attouch2013}) is     
  \vspace{-0.1in}$$u_i(x_i,z) = f(x_i,z_{\ne i}) + \kappa_i D_{\mathfrak h_i} (x_i,z_i),     
    $$
    where $\kappa_i$ is a scalar that can depend on $z$, and $D_{\mathfrak h_i} (x_i,z_i)$ is a Bregman divergence defined by 
    \vspace{-0.1in} \begin{equation}
        \label{def:Bregman}
        D_{\mathfrak h_i}(x_i,z_i)=\mathfrak h_i(x_i) - \mathfrak h_i(z_i) - \langle \nabla \mathfrak h_i(z_i),x_i - z_i \rangle,       
    \end{equation} 
     where $\mathfrak h_i:\mathbb R^{n_i}\to \mathbb R$ is a strongly convex function and can be adaptively chosen in the course of the update of $x_i$. If $\mathfrak h_i(x_i) = \frac12\langle x_i, Q_i x_i \rangle$, where   $Q_i$ is a positive definite matrix that could also depend on $z$,  then the proximal surrogate becomes the typical extended proximal surrogate       
    \begin{equation*}
      u_i(x_i,z) =f(x_i,z_{\ne i}) + \frac12\kappa_i \|x_i-z_i\|_{Q_i}^2.          \end{equation*}
  The approximation error 
  is $e_i(x_i,z)=\kappa_i D_{\mathfrak h_i} (x_i,z_i).$
    
    \item[(ii)] A quadratic surrogate (see e.g., \cite{Emilie2016,Ochs2019}) is defined as     
    $$u_i(x_i,z) = f(z) + \langle \nabla_{x_i} f(z), x_i-z_i \rangle + \frac{\kappa_i}{2} \|x_i-z_i\|_{Q_i}^2,  
    $$
    where $\kappa_i\geq 1$, $f$ is assumed to be twice differentiable, and $Q_i$ (which can depend on $z$) is a positive definite matrix such that $Q_i-\nabla_{x_i}^2 f (x_i,z_{\ne i})$ is also a positive definite matrix. The approximation error function is 
      \vspace{-0.1in}$$e_i(x_i,z)= f(z) - f(x_i,z_{\ne i}) + \langle \nabla_{x_i}f(z), x_i-z_i \rangle + \frac{\kappa_i}{2} \|x_i-z_i\|_{Q_i}^2.    
    $$ 
    If $Q_i=L_i I$, where $I$ is a identity matrix and $L_i$ is a positive number, then the  quadratic surrogate reduces to a Lipschitz gradient surrogate in the following.
    \item[(iii)] A Lipschitz gradient surrogate (see e.g.,  \cite{Xu2013}) is defined as
        $$u_i(x_i,z) = f(z) + \langle \nabla_{x_i} f(z), x_i-z_i \rangle + \frac{\kappa_i L_i}{2} \|x_i-z_i\|^2,    
    $$ 
    where $\kappa_i\geq 1$ and we assume $x_i \mapsto \nabla_{x_i} f(x_i,z_{\ne i})$ is $L_i$-Lipschitz continuous. Note that $L_i$ can depend on $z$.
    
    \item[(iv)] A Bregman surrogate (see e.g., \cite{HienNicolas_KLNMF,Hien2022,melo2017}) is defined as 
    \vspace{-0.1in} $$u_i(x_i,z) = f(z) + \langle \nabla_{x_i} f(z), x_i-z_i \rangle + \kappa_i L_i D_{\mathfrak h_i}(x_i,z_i),    
    $$ 
    where $\kappa_i\geq 1$, $\mathfrak h_i:\mathbb R^{n_i}\to \mathbb R$ is a strongly convex function, $D_{\mathfrak h_i}(x_i,z_i)$ is a Bregman divergence defined in \eqref{def:Bregman}, the function $x_i \mapsto f(x_i,z_{\ne i})$ is assumed to be $L_i$-relative smooth ($L_i$ can depend on $z$) to $\mathfrak h_i$, that is, the function $x_i\mapsto L_i \mathfrak h_i(x_i) -  f(x_i,z_{\ne i})$ is convex (see \cite{lu2018relatively}).      
    If $\mathfrak h_i(x_i)=\frac12\|x_i\|^2$ then the block Bregman surrogate reduces to the Lipschitz gradient surrogate.  
\end{enumerate}
\end{example}

 \section{Multiblock ADMM for solving Problem~\eqref{eq:model}} 
 \label{sec:mADMM}
 
We describe \mADMM for solving Problem~\eqref{eq:model} in Algorithm \ref{alg:mADMM-linear}.    
    \begin{algorithm}[tb]
\caption{mADMM for solving Problem~\eqref{eq:model}} 
\label{alg:mADMM-linear} 
\begin{algorithmic} 
\STATE \textbf{Notations.} For a given $k$ (the outer iteration index) and $i$ (the cyclic inner iteration index) we denote 
$x^{k,i}=(x^{k+1}_1, \ldots,x^{k+1}_{i},x^{k}_{i+1},\ldots,x^{k}_s)$ and let $x^{k+1} = x^{k,m}$.
\STATE Choose  $x^{0}$, $y^0$, $\omega^0$ and $\beta$ satisfies \eqref{para_condition_2}.
Set $k=0$.
\REPEAT
\FOR{ $i = 1,...,m$}
   \STATE  
   Update 
    \vspace{-0.1in} \begin{equation}
       \label{eq:x_iupdate-2}
       x^{k+1}_i\in \arg\min_{x_i} \{ u_i(x_i,x^{k,i-1},y^k,\omega^k) + g_i(x_i) \}. 
   \end{equation} 
   \ENDFOR
   \STATE Update $y$ as in~\eqref{eq:y_update}.
     \STATE Update $\omega$ as follows  
          \vspace{-0.1in} \begin{equation}
   \label{eq:omega_update-2}
   \omega^{k+1}=\omega^k  + \beta \big(\mathcal  \phi( x^{k+1}) + \mB y^{k+1}\big). 
   \end{equation} 
    \STATE Set $k\leftarrow k+1$.
\UNTIL{some stopping criteria is satisfied.}
\end{algorithmic}
\end{algorithm} 
We note that $x_i$, $i\in [m]$, $y$ and $\omega$ are the blocks of variables of $\varphi_\beta$ defined in \eqref{Lagrangian}.  

\textit{Update of block $x_i$.} We choose block surrogate functions of $\varphi_\beta$ with respect to $x_i$,  $i\in [m]$, such that they satisfy Condition~\ref{condition-C3} and one of the two conditions - Condition~\ref{condition-C1} or Condition~\ref{condition-C2}. 
\begin{condition}
\label{condition-C3}
     For $i\in [m]$, there exists an ``upper bound" error $ \bar{e}_i(x_i, z)$ such that the approximation error $e_{i}(x_i, z)=u_{i}(x_i, z)-\varphi_\beta (x_i, z_{\ne i})$ satisfies 
\begin{itemize}
\item[$\bullet$] $e_{i}(x_i, z) \leq \bar{e}_i(x_i, z)$ for all $(x_i, z)$, i.e.,  $e_{i}$ is upper bounded by $\bar e_i$, and
\item[$\bullet$]  we have $\bar{e}_i(z_i, z)=0$  and
$
    \nabla_{x_i} \bar e_{i}(  z_i,  z)=0$  for all $z$. 
\end{itemize}
\end{condition}

\begin{condition}
\label{condition-C1} 
The error $e_i(x_i,z)$ satisfies $e_i(x_i,z) \geq \eta_i(z,\beta) D_i(x_i,z_i)$ for all $x_i$, $z$, where $\eta_i(z,\beta)$ is a scalar that can depend on the values of $z$ and $\beta$. Here $D_i(\cdot,\cdot)$  satisfies 
\vspace{-0.1in} \begin{equation}
\label{condition-Di}
D_i(x_i^{k+1},x_i^{k}) \geq 0, \, \mbox{and if} \, D_i(x_i^{k+1},x_i^{k})\to 0 \,  \mbox{when} \,k\to\infty \, \mbox{then} \,\|x_i^{k+1}-x_i^k\|\to 0.
\end{equation} 
A simple example of $D_i(\cdot,\cdot)$ is $D_i(x_i^{k+1},x_i^{k})=\frac12\|x_i^{k+1}-x_i^k\|^2$.  
\end{condition}

\begin{condition}
\label{condition-C2} 
 The function $t(x_i,z)=u_i(x_i,z) + g_i (x_i)$ satisfies
$$ t(x_i,z) \geq t(x_i',z) + \langle \xi, x_i-x_i'\rangle + \eta_i(z,\beta) D_i(x_i',x_i)$$ for all $x_i,x_i',z$  and any $\xi \in \partial_{x_i} t(x_i',z)$. Here $D_i(\cdot,\cdot)$  satisfies the conditions in \eqref{condition-Di}. For the case $D_i(x_i',x_i)=1/2\|x_i'-x_i\|^2$, Condition \ref{condition-C2} reduces to the strong convexity of $x_i \mapsto t(x_i,z)$ with constant $\eta_i(z,\beta)$.   
\end{condition}

\begin{remark} 
\label{remark1} Note that all surrogates in Example \ref{example:surrogates} satisfy Condition \ref{condition-C3} with $\bar{e}_i(x_i, z)=e_{i}(x_i, z)$.

The surrogate in Example \ref{example:surrogates}(i) satisfies Condition \ref{condition-C1} with $\eta_i=\kappa_i$ and $D_i=D_{\mathfrak h_i}$. Consider the Bregman surrogate in  Example \ref{example:surrogates}(iv). Since  $x_i\mapsto L_i \mathfrak h_i(x_i) -  f(x_i,z_{\ne i})$ is convex, we have 
\vspace{-0.1in}$$ L_i \mathfrak h_i(x_i) -  f(x_i,z_{\ne i})- (L_i \mathfrak h_i(z_i) -  f(z_i,z_{\ne i}))-\langle L_i \nabla \mathfrak h_i(z_i)-\nabla_{x_i} f(z_i,z_{\ne i}),x_i-z_i\rangle \geq 0,
$$
which implies 
$L_i D_{\mathfrak h_i}(x_i,z_i) + f(z)+  \langle \nabla_{x_i} f(z), x_i-z_i \rangle- f(x_i,z_{\ne i})\geq 0. $
Therefore, if we take $\kappa_i>1$ then the approximation error of the Bregman surrogate satisfies Condition~\ref{condition-C1} with $\eta_i=(\kappa_i-1) L_i$ and $D_i=D_{\mathfrak h_i}$. Specifically, 
 \vspace{-0.1in} \[
e_i(x_i,z)= f(z) + \langle \nabla_{x_i} f(z), x_i-z_i \rangle + \kappa_i L_i D_{\mathfrak h_i}(x_i,z_i)-f(x_i,z_{\ne i})\geq (\kappa_i-1) L_i D_{\mathfrak h_i}(x_i,z_i).
 \] 
Similarly, we can show that the quadratic surrogate and the Lipschitz gradient surrogate in Example~\ref{example:surrogates} also satisfy Condition~\ref{condition-C1}.

Consider the Lipschitz gradient surrogate in Example~\ref{example:surrogates} (iii). If $g_i$ is convex then we see that $x_i\mapsto u_i(x_i,z) + g_i(x_i)$ is $\kappa_i L_i$-strongly convex. Hence, if $g_i$ is convex then Condition~\ref{condition-C2} is satisfied with $\eta_i(z,\beta)=L_i(z,\beta)$ and $D_i(x_i',x_i)=\frac12\|x_i'-x_i\|^2$. Similarly, we can prove that if $g_i$ is convex, then the quadratic surrogate and the Bregman surrogate in Example~\ref{example:surrogates} also satisfy Condition~\ref{condition-C2}. 
\end{remark}

\textit{Update of block $y$.}
As $h$ is assumed to be $L_h$-smooth, we use the following surrogate
\begin{equation}
\label{um-surrogate} u_{m+1}(y,\tilde x,\tilde y,\omega) = h(\tilde y)+  \langle \nabla h(\tilde y), y-\tilde y \rangle + \frac{ L_h}{2} \|y-\tilde y\|^2 +  \bar\varphi_\beta(\tilde x,y,\omega),
\end{equation} 
where $ \bar\varphi_\beta(x,y,\omega):= f(x) +   \langle \omega,\phi(x) + \mcalB y \rangle + \frac{\beta}{2} \| \phi(x) + \mcalB y\|^2$.
 The update of $y$ is
\vspace{-0.1in} \begin{equation}
\label{eq:y_update}
\begin{array}{ll}
     &  y^{k+1}\in\arg\min_{y}  u_{m+1}(y,x^{k+1},y^k,\omega^k)\\
       &= \arg\min\limits_{y}\big\{ \langle\nabla h(y^k) + \mB^*\omega^k ,y \rangle + \frac{\beta}{2} \| \phi( x^{k+1}) + \mB y \|^2+\frac{L_h}{2}\|y- y^k\|^2\big\}\\
     &=(\beta \mB^*\mB + L_h \mathbf I)^{-1}\big(L_h y^k - \nabla h(y^k)-\mB^*(\omega^k + \beta( \phi(x^{k+1})))\big).
       \end{array}
 \end{equation}  

\subsection{Subsequential convergence.} We establish subsequential convergence for 
 \mADMM\,in this section. Let us first prove some sufficient decreasing properties for the update of $x_i$ and $y$.
\begin{proposition}
\label{prop:nsdp-sufficient-xi}
(i) The update in \eqref{eq:x_iupdate-2} guarantees a sufficient decreasing:
\vspace{-0.1in} \begin{equation}
\label{NSDP}
\mL_{\beta}(x^{k,i},y^k,\omega^k) + \eta_i^k D_i( x^{k+1}_i, x_i^{k})
\leq
\mL_{\beta}(x^{k,i-1},y^k,\omega^k),
\end{equation} 
where $\eta_i^k=\eta_i(x^{k,i},y^k,\omega^k,\beta)$, which is the scalar in Condition \ref{condition-C1} or Condition \ref{condition-C2}. 

 (ii) The update in \eqref{eq:y_update} guarantees a sufficient decreasing:   
  \vspace{-0.1in} \begin{equation}
    \label{nsdp_y}
    \mL_{\beta}(x^{k+1},y^{k+1},\omega^k) + \frac{\delta}{2} \| y^{k+1} - y^k\|^2\leq \mL_{\beta}(x^{k+1},y^{k},\omega^k),   
\end{equation}
where  $\delta=L_h + \beta \lambda_{\min}(\mB^*\mB)$.
\end{proposition}
\proof{Proof.}
     We remark that, if we denote $x_{m+1}=y$, $g_{m+1}(y)=0$ and $D_{m+1}(y,y')=\frac12 \|y-y'\|^2$, then   
we can write the update of block $x_i$, $i\in[m]$,  and the update of block $y$ in \eqref{eq:y_update} in the following unified form, supposed the current iterate is  $ z=(\tilde x_1, \ldots,\tilde x_{m+1}, \omega)$,
  \begin{equation}
    \label{eq:updatex_i_rewrite}
x_i^+=\arg\min_{x_i} \{u_i(x_i, z) + g_i(x_i)\},
\end{equation}  
where $u_i$ is a surrogate function of $\varphi_\beta$ with respect to block $x_i$, $i\in [m+1]$.

(i) In the following, we use $\eta_i$ for $\eta_i(z,\beta)$. Suppose Condition \ref{condition-C1} holds. Then we have 
\vspace{-0.1in} \begin{equation}
   \label{temp:cal1}
   u_i(x_i^+,z)- \varphi_\beta(x_i^+,z_{\ne i}) \geq \eta_i D_i(x_i^+,z_i). 
\end{equation} 
On the other hand, from \eqref{eq:updatex_i_rewrite} we have
$  
u_i(x_i,z) + g_i(x_i) \geq u_i(x_i^+,z) + g_i(x_i^+)$ for all $x_i.$
Choose $x_i=z_i$, then combine with \eqref{temp:cal1}. We obtain
 \vspace{-0.1in} \begin{equation}
   \label{temp:cal3}
   u_i(z_i,z) + g_i(z_i) \geq u_i(x_i^+,z) + g_i(x_i^+) \geq \varphi_\beta(x_i^+,z_{\ne i}) + \eta_i D_i(x^+_i,z_i) + g_i(x_i^+).  
\end{equation} 
Furthermore,  $u_i(z_i,z)=\varphi_\beta(z)$. The result follows from \eqref{temp:cal3}.

Suppose Condition \ref{condition-C2} holds. Choosing $x_i'=x_i^+$ and $x_i=z_i$ in  Condition \ref{condition-C2}, we have 
  \begin{equation}
   \label{temp:cal4}
   u_i(z_i,z) + g_i(z_i) \geq t(x_i^+,z) + \langle \xi, z_i - x_i^+\rangle + \eta_i  D_i (x_i^+,z_i), \quad\forall \xi\in\partial_{x_i} t(x_i^+,z). 
\end{equation} 
From \eqref{eq:updatex_i_rewrite} we have $0\in \partial_{x_i} t(x_i^+,z)$. Hence, it follows from \eqref{temp:cal4} that
$$\varphi_\beta(z) +  g_i(z_i) \geq u_i(x_i^+,z) + g_i(x_i^+) + \eta_i  D_i (x_i^+,z_i).
$$ 
Moreover, we have $u_i(x_i^+,z) \geq \varphi_\beta(x_i^+,z_{\ne i})$. Hence the result follows.

(ii) Since $y\mapsto u_{m+1}(y,\tilde x,\tilde y,\omega)$ defined in \eqref{um-surrogate} is  $(L_h+ \beta \lambda_{\min}(\mB^*\mB))$-strongly convex, similarly to the proof of Part (i) for the case that Condition~\ref{condition-C2} is satisfied, we can prove that \eqref{nsdp_y} is satisfied.
\QED

Denote $\Delta x^k_i= x^{k}_i -x^{k-1}_i$, $\Delta y^k=y^{k} -y^{k-1}$, $\Delta \omega^k =\omega^{k} -\omega^{k-1}$.

\begin{proposition}
\label{prop:recursive1}
 The values of $\eta_i^k$ and $\delta$ are given in Proposition~\ref{prop:nsdp-sufficient-xi}. Let $\hat\delta=2L_h$. We have
\vspace{-0.1in} \begin{equation}
    \label{nsdp-L}
    \begin{array}{ll}
   &\mL_{\beta}(x^{k+1},y^{k+1},\omega^{k+1}) + \sum\limits_{i=1}^m\eta^k_i D_i(x^{k+1}_i,x^k_i) 
   + \frac{\delta}{2} \| \Delta y^{k+1}\|^2 \\  
   & \leq \mL_{\beta}(x^{k},y^{k},\omega^{k})  + \frac{1}{\beta \sigma_{\mcalB}}  \| \mcalB^* \Delta \omega^{k+1} \|^2,
   \end{array} 
\end{equation} 
and 
 \begin{equation}
\label{ieq5}
 \|\mcalB^* \Delta\omega^{k+1}\|^2 \leq 3 \left( (L_h^2 +\hat \delta^2) \| \Delta y^{k+1}\|^2 +  \hat \delta^2\| \Delta y^{k}\|^2  \right).
\end{equation} 
\end{proposition}

\proof{Proof.}
    Summing \eqref{NSDP} from $i=1$ to $m$ we obtain
 \begin{equation}
    \label{ie1}
    \mL_{\beta}(x^{k+1},y^k,\omega^k) + \sum_{i=1}^m\eta^k_i D_i(x^{k+1}_i,x^k_i)
\leq
\mL_{\beta}(x^{k},y^k,\omega^k). 
\end{equation} 
On the other hand, 
 \begin{equation}
\label{ie2}
\begin{array}{ll}
\mL_{\beta}(x^{k+1},y^{k+1},\omega^{k+1})  &= \mL_{\beta}(x^{k+1},y^{k+1},\omega^k) + \left\langle \omega^{k+1}-\omega^k, \phi(x^{k+1})+\mcalB y^{k+1} \right\rangle\\
&= \mL_{\beta}(x^{k+1},y^{k+1},\omega^k)
+ \frac{1}{\beta}\| \Delta \omega^{k+1}\|^2 \\
& \leq \mL_{\beta}(x^{k+1},y^{k+1},\omega^k)
+ \frac{1}{\beta \sigma_{\mcalB}} \| \mcalB^*\Delta \omega^{k+1}\|^2. 
\end{array}
\end{equation} 
Therefore, Inequality \eqref{nsdp-L} follows from~\eqref{ie1}, \eqref{nsdp_y}, and \eqref{ie2}.

By definition of the approximation error, we have 
$$\nabla_y e_{m+1}(y^{k+1}, x^{k+1}, y^k,\omega^k)=\nabla_y u_{m+1}(y^{k+1}, x^{k+1},y^k,\omega^k) - \nabla_y \varphi_\beta ( x^{k+1},y^{k+1},\omega^k),
$$ 
where $u_{m+1}$ is defined in \eqref{um-surrogate}.
Furthermore, from \eqref{eq:y_update}, 
 $\nabla_y u_{m+1}(y^{k+1}, x^{k+1}, y^k,\omega^k)=0.$ Hence, 
 $-\nabla_y e_{m+1}(y^{k+1}, x^{k+1}, y^k,\omega^k) = \nabla_y \varphi_\beta ( x^{k+1},y^{k+1},\omega^k),
 $ 
which implies that 
\begin{equation}
 \label{eq:connect_w}
 \begin{split}
-\nabla_y e_{m+1}(y^{k+1}, x^{k+1}, y^k,\omega^k) &= \nabla h(y^{k+1}) + \mcalB^*\big( \omega^k + \beta  (\phi(x^{k+1}) + \mcalB y^{k+1}) \big)\\
&=\nabla h(y^{k+1}) +  \mcalB^* \omega^{k+1}.
\end{split}
\end{equation}
Hence, we have 
 \begin{equation}
\label{ie4}
\begin{array}{ll}
\|\mcalB^*\Delta \omega^{k+1}\|&= \big \|-\nabla_y e_{m+1}(y^{k+1}, x^{k+1}, y^k,\omega^k) -\nabla h(y^{k+1})\\
&\qquad\qquad\qquad+ \nabla_y e_{m+1}(y^{k}, x^{k}, y^{k-1},\omega^{k-1})+\nabla h(y^{k})\big \|.
\end{array} 
\end{equation} 
On the other hand, $\nabla_y e_{m+1}(y,\tilde x,\tilde y,\omega)=\nabla h(\tilde y) - \nabla h(y)+L_h (y-\tilde y)$. Hence, 
  \vspace{-0.1in} \begin{equation}
\label{eq:errorh-2} \|\nabla_y e_{m+1}(y^{k+1},x^{k+1},y^k,\omega^k) \| \leq \hat\delta \|  y^{k+1}-y^k\|.
      \end{equation} 
So from~\eqref{ie4}  we obtain Inequality \eqref{ieq5}. 
\QED
\begin{proposition}
\label{prop:delta_converge}
Let $\delta$ and $ \eta^k_i$ be defined in Proposition~\ref{prop:nsdp-sufficient-xi} and $\hat\delta=2L_h$. Suppose  
     \vspace{-0.1in} \begin{equation}
     \label{para_condition_2}
     \beta(L_h + \beta \lambda_{\min}(\mB^*\mB)) \geq \frac{6L_h^2(5+4\tilde \delta)}{ \sigma_{\mathcal B}} \quad \mbox{and} \quad \eta^k_i\geq\underline{\eta}_i
     \end{equation}  
for some constants $\underline{\eta}_i>0$ and $\tilde\delta>1$.
Denote $\mL^k=\mL_{\beta}(x^{k},y^{k},\omega^{k})$. 

(A) For $k>0$, we have 
\vspace{-0.1in} \begin{equation}
\label{recursive-2}
\mL^{k+1} + \sum_{i=1}^m\eta^k_i D_i(x^{k+1}_i,x^k_i) + \frac{3\tilde\delta \hat\delta^2}{ \beta\sigma_{\mcalB}} \| \Delta y^{k+1}\|^2 
\leq \mL^k+   \frac{3 \hat\delta^2}{ \beta\sigma_{\mcalB}} \| \Delta y^{k}\|^2. 
\end{equation} 

(B) For $K>0$, we have
\vspace{-0.1in} \begin{equation}
\label{main_recursive}
\begin{array}{ll}
 & \mL^{K+1}+\frac{3\tilde\delta \hat\delta^2}{ \beta\sigma_{\mcalB}}   \| \Delta y^{K+1}\|^2 + \sum_{k=1}^K \sum_{i=1}^m\underline{\eta}_i D_i(x^{k+1}_i,x^k_i) \\
   &
  \qquad +  \frac{3(\tilde\delta -1)\hat\delta^2}{ \beta\sigma_{\mcalB}}  \sum_{k=1}^{K-1}  \| \Delta y^{k+1}\|^2 \leq  \mL^{1} +  \frac{3\hat\delta^2}{ \beta\sigma_{\mcalB}} \| \Delta y^{1}\|^2.  
\end{array}
\end{equation}  

(C) The sequences  $\{\Delta x^k\}$, $\{\Delta y^k\}$ and $\{\Delta \omega^k\}$ converge to 0.
\end{proposition}
\proof{Proof.}
(A)     Combining~\eqref{nsdp-L} with~\eqref{ieq5} gives us   
 \begin{equation*} 
\begin{array}{ll}
    &\mL^{k+1} + \sum_{i=1}^m\eta^k_i D_i(x^{k+1}_i,x^k_i) 
   + \frac{\delta}{2} \| \Delta y^{k+1}\|^2 \\
   &\quad\leq \mL^{k}  +  \frac{3}{\beta \sigma_{\mcalB}} \big( (L_h^2+\hat\delta^2)  \| \Delta y^{k+1}\|^2 +  \hat\delta^2\| \Delta y^{k}\|^2  \big).   
\end{array}
\end{equation*} 
This implies the following recursive inequality
\vspace{-0.1in} \begin{equation} 
\label{recuresive-1}
   \mL^{k+1} + \sum_{i=1}^m\eta^k_i D_i(x^{k+1}_i,x^k_i) 
   + \Big(\frac{\delta}{2} -\frac{3}{\beta \sigma_{\mcalB}}(L_h^2+\hat\delta^2 )\Big)\| \Delta y^{k+1}\|^2\leq \mL^{k}  + \frac{3}{\beta \sigma_{\mcalB}} \hat\delta^2\| \Delta y^{k}\|^2.
\end{equation}  
The condition in \eqref{para_condition_2} implies that 
$\frac{\delta}{2}-\frac{3}{\beta \sigma_{\mcalB}} \big(L_h^2+ \hat\delta^2   \big)\geq \frac{3\tilde\delta \hat\delta^2}{ \beta\sigma_{\mcalB}}. 
$ 
Hence Inequality~\eqref{recuresive-1} implies \eqref{recursive-2}.

 (B) Summing \eqref{recursive-2} from $k=1$ to $K$ and using the conditions $\eta_i^k\geq\underline{\eta}_i$ we get \eqref{main_recursive}. 
 
 (C) We now use the technique in \cite[Proposition 2.9]{melo2017} to prove that
\vspace{-0.1in} \begin{equation}
\label{eq:nu}
\hat\mL^k:=\mL^k+  \frac{3\tilde\delta \hat\delta^2}{ \beta\sigma_{\mcalB}}\| \Delta y^{k}\|^2\geq \nu, \quad\mbox{for all}  \quad k\geq 1,
\end{equation} 
where $\nu$ is a lower bound of $F(x^k)+h(y^k)$. From \eqref{recursive-2}, $\{\hat\mL^k\}_{k\geq 1}$ is non-increasing. Suppose there exists $k_1 \geq 1 $ s.t $\hat \mL^{k_1}<\nu$.  Hence,  $\hat \mL^k<\nu$ for all $k\geq k_1$. We have
$$
\sum_{k=1}^K(\hat{\mL}^k-\nu) \leq \sum_{k=1}^{k_1} (\hat{\mL}^k-\nu)  + (K-k_1)(\hat{\mL}^{k_1}-\nu).  
$$
Hence, $\lim_{K\to\infty}\sum_{k=1}^{K}(\hat{\mL}^k-\nu)=-\infty$. On the other hand,  we have 
$$
\hat{\mL}^k\geq \mL^k\geq F(x^k)+h(y^k)+\frac{1}{\beta}\langle \omega^k,\omega^k-\omega^{k-1}\rangle\geq \nu + \frac{1}{2\beta}(\|\omega^k\|^2-\|\omega^{k-1}\|^2). 
$$
Thus,$
\sum_{k=1}^K(\hat{\mL}^k-\nu) \geq \sum_{k=1}^K \frac{1}{2\beta}(\|\omega^k\|^2-\|\omega^{k-1}\|^2)=\frac{1}{2\beta}(\|\omega^K\|^2-\|\omega^{0}\|^2)\geq -\frac{1}{2\beta}\|\omega^{0}\|^2. 
$
This gives a contradiction. Therefore,  we get \eqref{eq:nu}.

Inequality \eqref{main_recursive} together with the lower boundedness of $\{\hat\mL\}_{k\geq 0}$ we have  $\sum_{k=1}^\infty D_i(x^{k+1}_i,x^k_i)<+\infty$, $\sum_{k=1}^\infty\| \Delta y^{k+1}\| <+\infty $. This implies that $\Delta x^{k}$ and $\Delta y^{k}$ converge to 0. 
On the other hand, it follows from \eqref{ieq5} that $\sum_{k=0}^\infty \|\mcalB^* \Delta\omega^{k+1}\|^2<+\infty$, leading to the convergence of $\{\Delta\omega^k\}$ to 0. 
\QED
\begin{theorem}[Subsequential convergence]
\label{thm:subsequential}
Suppose the parameters are chosen such that the conditions of Proposition~\ref{prop:delta_converge} are satisfied. If there exists a subsequence $(x^{k_n},y^{k_n},\omega^{k_n})$ converging to $(x^*,y^*,\omega^*)$ then $(x^*,y^*,\omega^*)$ is a critical point of $\mL_{\beta}(x,y,\omega)$. 
\end{theorem}
\proof{Proof.}
    See Apendix \ref{thrm-sub}.
\QED

In the following proposition,  by  extending \cite[Lemma 6]{Wang2019}, we provide a sufficient condition such that the generated sequence of \mADMM\, is bounded.
\begin{proposition}
\label{prop:bounded-sequence}
If $ran \, \phi(x)  \subseteq Im(\mB)$, $\lambda_{\min}(\mB^* \mB)>0 $ and $F(x) + h(y)$ is coercive over the feasible set $\{(x,y): \phi(x) + \mB y = 0\}$ then  $\{(x^k,y^k,\omega^k)\}_{k\geq0}$ generated by Algorithm~\ref{alg:mADMM-linear} is bounded.
\end{proposition}
\proof{Proof.}
    See Appendix \ref{bounded-sequence}.
\QED
\subsection{Iteration complexity.}
We now establish the iteration complexity to obtain an $\varepsilon$-stationary point, see Definition~\ref{def:epsilon_opt}. To this end, we need the following additional assumption. 

\begin{assumption}
\label{assump:u_i}
\begin{itemize}
\item For any $z$ and $x_i\in {\rm dom} (g_i)$, we have 
\begin{equation}
\label{assume:partial}
\begin{split}
\partial_{x_i} \big(u_i(x_i,z) + g_i(x_i) \big) &= \partial_{x_i} u_i(x_i,z) + \partial g_i(x_i),
\partial_{x_i} (f(x) + g_i(x_i))\\
&= \partial_{x_i} f(x) + \partial g_i(x_i). 
\end{split}
\end{equation}  
\item  For any $S^k_i \in\partial_{x_i} u_i(x_i^{k+1},x^{k,i-1},y^k,\omega^k)$ there exists $\bar S^k_i\in \partial_{x_i} \varphi_\beta(x^{k+1},y^k,\omega^k)$ that
\vspace{-0.1in} \begin{equation}
\label{eq:l1}
\|S_i^k - \bar S_i^k\|\leq \bar L_i \|x^{k+1}-x^{k,i-1}\|,
\end{equation} 
where $\bar L_i$ is some nonnegative constant (we remark that the constant $\bar L_i$ does not involve in how to choose the parameters in our framework; its existence is for the convergence proof). 
\end{itemize}
\end{assumption}
Condition \eqref{assume:partial} says that $x_i\mapsto u_i(x_i,z)$, $x_i\mapsto f(x)$ and $g_i$ follow the sum rule for the limiting subgradients. See  \cite[Corollary 10.9]{RockWets98} for a sufficient condition. 
Note that if $x_i\mapsto u_i(x_i,z)$ and $x_i\mapsto f(x)$ are continuously differentiable then \eqref{assume:partial} is satisfied.
Regarding~\eqref{eq:l1},  if we assume  $(x^k,y^k,\omega^k)$ is bounded then~\eqref{eq:l1} is satisfied when $\nabla_{x_i}u_i(x_i,x,y,\omega)=\nabla_{x_i}\varphi_\beta (x,y,\omega)$  and $u_i$ is twice continuously differentiable.  Indeed, consider the bounded set containing $(x^k,y^k,\omega^k)$, then \vspace{-0.1in}
\begin{equation*}
\begin{array}{ll}
\|S_i^k - \bar S_i^k\|&=\|\nabla_{x_i}u_i(x^{k+1}_i,x^{k,i-1},y^{k},\omega^{k})-\nabla_{x_i}\varphi_\beta (x^{k+1},y^{k},\omega^{k})\|\\
&=\|\nabla_{x_i}u_i(x^{k+1}_i,x^{k,i-1},y^{k},\omega^{k})-\nabla_{x_i}u_i(x_i^{k+1},x^{k+1},y^{k},\omega^{k})\|\\
&\leq \bar L_i \|x^{k+1}-x^{k,i-1}\|.
\end{array}
\end{equation*}
The surrogate functions given in Example~\ref{example:surrogates} satisfy~\eqref{eq:l1} when $f$ is twice continuously differentiable. We give  an example when $f$ is nonsmooth and~\eqref{eq:l1} is still satisfied in Appendix \ref{appendix:example}.

\begin{proposition}
\label{prop:iteration_complexity}
Denote 
$$R^k_y=\|\nabla h(y^k)+\mcalB^* \omega^k\|,  R_i^{k} = \|\chi^k_i + \nabla_{x_i} \phi(x^k)\omega^k\|,  R_c^k=\|\phi(x^k)+\mcalB y^k\|,
$$
where $\chi^{k} _i\in\partial_{x_i} F(x^{k})$. Assume  Assumption~\ref{assump:u_i} hold, the conditions of Proposition~\ref{prop:delta_converge} are satisfied, and $D_i(x_i^{k+1},x_i^{k})=\frac12\|x^{k+1}_i-x_i^{k}\|^2$ in \eqref{NSDP}. Suppose $\{(x^k,y^k)\}_{k\geq 0}$ is bounded. For all $k> 0$, there exist $\chi^{k} _i\in\partial_{x_i} F(x^{k})$ such that
\vspace{-0.1in} \begin{equation}
\label{R_i_rate}
R^k_y= O(\|\Delta y^{k}\|), R^{k}_i  = O(\|\Delta x^{k}\|+\|\Delta y^{k}\|), R^{k}_c = O(\|\Delta y^{k}\|+\|\Delta y^{k-1}\|) 
\end{equation} 
and
 \begin{equation}
\label{min_rate}
\min_{1\leq j\leq k-1} \Big\{ \frac12\tilde\eta \| \Delta x^{j+1}\|^2 
   + \frac{3(\tilde\delta-1)\hat\delta^2}{2\beta \sigma_{\mcalB}}\big( \| \Delta y^{j+1}\|^2 + \| \Delta y^{j}\|^2 \big)\Big\} =O\big(\frac{1}{k}\big),   
\end{equation} 
where $\tilde \eta=\min_{i\in [m]}\underline{\eta}_i$.
\end{proposition}
\proof{Proof.}
     From \eqref{eq:connect_w} and \eqref{eq:errorh-2} we have 
  \vspace{-0.1in} \begin{equation}
\label{forrate1}
 R^k_y =\|\nabla_y e_{m+1}(y^k,x^k,y^{k-1},\omega^{k-1}) \|\leq \hat{\delta}\|\Delta y^k\|.
 \end{equation} 
Writing the optimality condition for \eqref{eq:x_iupdate-2} we get
\vspace{-0.1in}$$0\in \partial_{x_i} \big( u_i(x_i^{k+1},x^{k,i-1},y^k,\omega^k)+g_i(x_i^{k+1})\big)=\partial_{x_i} u_i(x_i^{k+1},x^{k,i-1},y^k,\omega^k) + \partial_{x_i} g_i(x_i^{k+1}). 
$$
In other words, there exist $S^k_i \in\partial_{x_i} u_i(x_i^{k+1},x^{k,i-1},y^k,\omega^k)$ and $T_i^{k+1}\in \partial_{x_i} g_i(x_i^{k+1})$ such that $ S^k_i + T_i^{k+1}=0$. Furthermore, it follows from \eqref{eq:l1} that there exists $\bar S^k_i\in \partial_{x_i} \varphi_\beta(x^{k+1},y^k,\omega^k)$ such that $\|S_i^k - \bar S_i^k\|\leq \bar L_i \|x^{k+1}-x^{k,i-1}\|$. Hence, we have 
 \begin{equation} 
\label{tiski}
\|T_i^{k+1} + \bar S^k_i\| = \|\bar S^k_i - S_i^k \| \leq \bar L_i \|x^{k+1}-x^{k,i-1}\| \leq  \bar L_i \|\Delta x^{k+1}\|.
\end{equation} 
Note that 
\[
\begin{array}{ll}
 \partial_{x_i} \varphi_{\beta}(x^{k+1},y^k,\omega^k)& = \partial_{x_i} f(x^{k+1})+ \nabla\phi(x^{k+1})\Big(\omega^k + \beta(\phi(x^{k+1})+\mcalB y^k) \Big)  \\
 &= \partial_{x_i} f(x^{k+1})+ \nabla\phi(x^{k+1})\omega^{k+1} +\beta \nabla\phi(x^{k+1})(\mcalB y^k-\mcalB y^{k+1}).
 \end{array}
  \] 
Therefore, $\exists \xi_i^{k+1}  \in\partial_{x_i} f(x^{k+1})$ s.t
$\xi_i^{k+1} + \nabla\phi(x^{k+1})\omega^{k+1} +\beta \nabla\phi(x^{k+1})(\mcalB y^k-\mcalB y^{k+1}) = \bar S^k_i.
$
 Moreover, it follows from \eqref{assume:partial} that $\chi_i^{k+1}=\xi_i^{k+1} +  T^{k+1}_i \in \partial_{x_i} F(x^{k+1})$. We then obtain
 \begin{equation}
\label{forrate2}
\begin{array}{ll}
R_i^{k+1} &= \|\xi_i^{k+1}  + \nabla\phi(x^{k+1})\omega^{k+1} + T^{k+1}_i \|\\
&=\|\bar S_i^k - \beta\nabla \phi(x^{k+1})(\mcalB(y^k-y^{k+1}))+ T^{k+1}_i\|\\
&\leq \bar L_i \|\Delta x^{k+1}\| + \beta\mathbf M_\phi \|\mcalB\|\|\Delta y^{k+1}\|,
 \end{array}
\end{equation}
 where in the last inequality we used: \eqref{tiski}, the fact that $ \|\nabla\phi(x^{k+1})\| \leq\mathbf M_\phi$ for some constant $\mathbf M_\phi$ (since $x^k$ is bounded),  $\phi$ is continuously differentiable. 
 
 On the other hand,  we have 
$$
 R^{k+1}_c =\|\phi(x^{k+1}) + \mcalB y^{k+1}\|=\frac{1}{\beta}\|\Delta \omega^{k+1}\| \leq \frac{1}{\beta\sigma_\mcalB} \| \mcalB^*  \Delta \omega^{k+1}\| 
.$$ Together with  \eqref{ieq5}, \eqref{forrate1}, and \eqref{forrate2},  we obtain \eqref{R_i_rate}.
 
From \eqref{main_recursive} and \eqref{eq:nu} we have  
 \begin{equation}
\label{eq:nu_eta_D}
  \nu + \sum_{k=1}^K \sum_{i=1}^m\underline{\eta}_i D_i(x^{k+1}_i,x^k_i) 
   +  \frac{3(\tilde\delta-1) \hat\delta^2}{\beta \sigma_\mcalB} \sum_{k=1}^{K-1}  \| \Delta y^{k+1}\|^2 \leq  \mL^{1} +  \frac{3 \hat\delta^2}{\beta \sigma_{\mcalB}} \| \Delta y^{1}\|^2. 
\end{equation} 
Note that $ 2 \sum_{k=1}^{K-1}  \| \Delta y^{k+1}\|^2 \geq \sum_{k=1}^{K-1} \big( \| \Delta y^{k+1}\|^2 +  \| \Delta y^{k}\|^2\big) -  \| \Delta y^{1}\|^2$. 
Therefore from \eqref{eq:nu_eta_D}, we have
$$
\min_{1\leq k\leq K-1} \Big\{\frac12\tilde{\eta}  \| \Delta x^{k+1}\|^2 
   + \frac{3(\tilde\delta-1)\hat\delta^2}{2\beta \sigma_{\mcalB}}  \big( \| \Delta y^{k+1}\|^2 +  \| \Delta y^{k}\|^2\big)\Big\} \leq \frac{1}{K-1}C, 
   $$ 
   where $C=( \mL^{1} + \frac{3 \hat\delta^2}{\beta \sigma_{\mcalB}}\| \Delta y^{1}\|^2-\nu + \frac{3(\tilde\delta-1)\hat\delta^2}{2\beta \sigma_{\mcalB}}\| \Delta y^{1}\|^2)$.
   The result in \eqref{min_rate} follows then.
\QED

Proposition \eqref{prop:iteration_complexity} shows that  
$\min_{0\leq j\leq k-1} R^j_y + \sum_{i=1}^m R^{j}_i + R^j_c = O\big(\frac{1}{\sqrt{k}}\big).
$
Hence, after at most $O\big( \frac{1}{\varepsilon^2}\big)$ iterations we obtain an $\varepsilon$-stationary point of~\eqref{eq:model}.

\subsection{Global convergence under K{\L} property and convergence rate.}

The following theorem presents the global convergence of mADMM under the K{\L} property (see \cite{Bolte2014} for its definition).  
\begin{theorem}
\label{thrm:global}
Assume Assumption~\ref{assump:u_i} hold, the conditions of Proposition~\ref{prop:delta_converge} are satisfied, and $D_i(x_i^{k+1},x_i^{k})=\frac12\|x^{k+1}_i-x_i^{k}\|^2$ in \eqref{NSDP}. Furthermore, we assume that the generated sequence is bounded and the following Lyapunov function 
 \vspace{-0.1in} \begin{equation}
\label{eq:lyapunov}
\tilde{\mL}_\beta (x,y,\omega,\tilde y) = \mL_\beta(x,y,\omega) +   \frac{3 \hat \delta^2}{\beta \sigma_{\mcalB}} \|y - \tilde y\|^2
\end{equation} 
has the K{\L}~property with constant $\sigma_{\tilde \mL}$, then $\{z^{k}\}_{k\geq0} $, where $z^k=(x^k,y^k,\omega^k)$, converges to a critical point of $\mL_{\beta}$. 
Moreover, if $\sigma_{\tilde L}=0$ then \mADMM\, converges after a finite number of steps; if $\sigma_{\tilde L} \in (0,1/2]$ then there exists $k_1\geq 1$, $w_1>0$ and $w_2\in [0,1)$ such that $ \|z^k-z^*\| \leq w_1 w_2^k$ for all $k\geq k_1$; and if $ \sigma_{\tilde L} \in (1/2,1]$ then there exists $k_1\geq 1$ and $w_1>0$ such that  $ \|z^k-z^*\| \leq w_1 k^{-(1-\sigma_{\tilde L})/(\sigma_{\tilde L} -1)}$  for all $k\geq k_1$. 
\end{theorem}

\proof{Proof.}
See Appendix~\ref{proof-global}
\QED

 \begin{remark}
The assumption $ \sigma_{\mB}=\lambda_{\min}(\mB \mB^*)>0$ can be relaxed to
$ ran \, \phi(x)  \subseteq Im(\mB).
$  
Then the constant $\sigma_{\mB}$ in \eqref{para_condition_2} would be replaced by \emph{the smallest positive eigenvalue} of $\mB^* \mB$, which is denoted by $\lambda^+_{\min} ( \mB^* \mB)$. It is because when $ran \, \phi(x)  \subseteq Im(\mB)$, we can derive from  \eqref{eq:omega_update-2} that $\Delta \omega^{k+1}\in  Im(\mB)$. This implies $ \|\Delta \omega^{k+1}\|^2 \leq \lambda^+_{\min} ( \mB^* \mB) \| \mB^*\Delta \omega^{k+1}\|^2 $, which is the cornerstone to prove Proposition \ref{prop:recursive1} (see Inequality \eqref{ie2}) and to derive $\Delta \omega^k\to 0$  from $\mB^* \Delta \omega^k\to 0$ in Proposition \ref{prop:delta_converge} (C).
 \end{remark}

 \section{Preliminary numerical results}
 \label{sec:numerical}


  We consider the following composite optimization problem 
   \vspace{-0.1in} \begin{equation}
 \label{regularized_nlq}
 \min_{x\in \mathbb R^n}  g(x) + h(\phi(x)),
 \end{equation} 
 where $g$ is a lower semicontinuous function, $h$ is a differentiable function with $L_h$-Lipschitz continuous gradient, and $\phi(x)=(\phi_1(x),\ldots,\phi_q(x))$ is a nonlinear mapping. 
Problem~\eqref{regularized_nlq}, which includes regularized nonlinear least square problem as a special case (that is, when $h(\phi(x))=\frac12 \| \phi(x)\|^2$), has appeared frequently in machine learning and statistics. For examples,  Problem~\eqref{regularized_nlq} covers the nonlinear regression problem \cite{Dutter1981} ($h$ plays the role of a loss function, $\phi$ represents the model to train and $g$ is a regularizer),  the risk parity portfolio selection problem \cite{Maillard2010}, the robust phase retrieval problem \cite{Duchi2018}, and  the PDE-constrained inverse problem  \cite{Roosta2014}.
 Problem~\eqref{regularized_nlq} is rewritten in the form of \eqref{eq:model}  as follows.
 \vspace{-0.1in} \begin{equation}
 \label{nlq_2}
   \min_{x\in \mathbb R^n, y\in\mathbb R^q}  g(x) + h(y)\quad \mbox{such that}\quad   \phi(x) - y = 0.
 \end{equation} 
To illustrate the effect of using block surrogate functions in \mADMM, let us consider the following specific example of \eqref{regularized_nlq}
 \begin{equation}
\label{svm_quad}
\min_{x=(x_1,x_2,x_3)} \frac{1}{q}\sum_{i=1}^q\log(1+ e^{-\mathbf b_i(\langle \mathbf a_i,x_1 \rangle^2 + \langle \mathbf a_i, x_2\rangle  + x_3)}) + \lambda\sum_{i=1}^2 \|x_i\|_1,
\end{equation} 
that is, $h(y)=\frac{1}{q}\sum_{i=1}^q\log(1+ e^{-\mathbf b_i y_i})$, $g(x)=\lambda_1 \|x_1\|_1 + \lambda_2 \|x_2\|_1$ and $ \phi_i(x_1,x_2,x_3)=\langle \mathbf a_i,x_1 \rangle^2 + \langle \mathbf a_i, x_2\rangle  + x_3$, where $\lambda_1$ and $\lambda_2$ are regularizer parameters, $\mathbf a_i\in\mathbb R^d$ and $\mathbf b_i \in\{-1,1\}$ are input data. Problem~\eqref{svm_quad} has a form of a support vector machine problem \cite{Kecman2005} that minimizes a logistic loss function with a nonlinear classifier and an $l_1$ regularization. 

\subsection{Applying mADMM to solve Problem~\eqref{svm_quad}.} 
The augmented Lagrangian for \eqref{nlq_2} is \vspace{-0.1in}
$$\mL_\beta(x,y,\omega)= g(x) +  h(y)+ \sum_{i=1}^q  \omega_i( \phi_i(x) - y_i ) + \frac{\beta}{2}\sum_{i=1}^q(\phi_i(x)-y_i)^2. 
$$ 

\textit{Update of block $x_1$.} Fix $\omega$ and $y$, and denote $\mathfrak f_i(x)= \omega_i( \phi_i(x) - y_i ) + \frac{\beta}{2} (\phi_i(x)-y_i)^2$.  We have 
\vspace{-0.1in}
$$
\begin{array}{ll} 
 \nabla_{x_1} \mathfrak f_i(x)&=2\omega_i \langle \mathbf a_i,x_1 \rangle \mathbf  a_i + 2\beta (\phi_i(x)-y_i) \langle \mathbf a_i,x_1 \rangle \mathbf a_i, \\
 \nabla^2_{x_1} \mathfrak f_i(x) & = 2\omega_i \mathbf  a_i \mathbf a_i^\top + 2 \beta (\phi_i(x)-y_i) \mathbf a_i \mathbf a_i^\top +  4 \beta   (\langle \mathbf a_i,x_1 \rangle)^2  \mathbf a_i  \mathbf a_i^\top.
 \end{array} 
 $$
This implies that

$$
\begin{array}{ll}
\|\nabla^2_{x_1} \mathfrak f_i(x) \| &\leq \| 2(\omega_i-\beta y_i  )\mathbf a_i \mathbf a_i^\top\| + 2\beta\|\mathbf a_i\mathbf  a_i^\top \| |\phi_i(x) + 2(\langle \mathbf a_i,x_1 \rangle)^2|\\
&\leq 2|\omega_i-\beta y_i  | \| \mathbf a_i\|^2 +  2\beta\|\mathbf a_i \|^2 (|\langle \mathbf a_i,x_2\rangle + x_3| + 3\|\mathbf a_i\|^2 \|x_1\|^2)\\
&= 2\big( |\omega_i-\beta y_i | + \beta  |\langle \mathbf a_i,x_2\rangle + x_3|\big) \|\mathbf  a_i\|^2 + 6\beta \|\mathbf  a_i\|^4 \|x_1\|^2\\
&\leq 2\|\mathbf  a_i\|^2 \max \Big\{\big( |\omega_i-\beta y_i | + \beta  |\langle \mathbf a_i,x_2\rangle + x_3|\big) , 3\beta \|\mathbf  a_i\|^2 \Big\} (1+\|x_1\|^2 ) . 
 \end{array} 
 $$
 It follows from~\cite[Proposition 2.1]{lu2018relatively} that $x_1\mapsto \mathfrak f_i(x)$ is $l_i$-relative smooth to 
 \vspace{-0.1in} \begin{equation}
 \label{h1}
 \mathfrak h(x_1)=\frac14 \|x_1\|_2^4 + \frac12 \|x_1\|_2^2 
  \end{equation}  
    with $l_i= 2\|\mathbf  a_i\|^2 \max \Big\{\big( |\omega_i-\beta y_i | + \beta  |\langle \mathbf a_i,x_2\rangle + x_3|\big) , 3\beta \| \mathbf a_i\|^2 \Big\}.$
Hence $x_1\mapsto \varphi_\beta(x,y,\omega)$ is $\sum_{i=1}^q l_i $-relative smooth to $\mathfrak h_i$ ($\varphi_\beta$ is defined in \eqref{Lagrangian}).  Denote $\mathfrak l_{1}=\sum_{i=1}^q l_i$.  We use the Bregman surrogate: 
\vspace{-0.1in} \begin{equation}
\label{eq:updatex_k}
\begin{array}{ll}
&x_1^{k+1}\in\argmin_{x_1} g_1(x_1) + \langle \nabla_{x_1}\varphi_\beta(x_1,x_2^k,x_3^k,y^k,\omega^k),x_1-x_1^k \rangle+ \mathfrak l_1 D_{\mathfrak h}(x_1,x_1^k)\\
&=\argmin_{x_1}\lambda_1 \|x_1\|_1 + \langle \nabla_{x_1}\varphi_\beta(x_1,x_2^k,x_3^k,y^k,\omega^k)-\mathfrak l_1\nabla \mathfrak h(x_1^k),x_1-x_1^k \rangle+ \mathfrak l_1 \mathfrak h(x_1)\\
&=\argmin_{x_1} \lambda_1 \|x_1\|_1 + \langle \tilde  c ,x_1\rangle + \mathfrak l_1 \mathfrak h(x_1), 
\end{array} 
\end{equation} 
where $\tilde c=\nabla_{x_1}\varphi_\beta(x_1,x_2^k,x_3^k,y^k,\omega^k)-\mathfrak l_1\nabla \mathfrak h(x_1^k)$.   
\begin{lemma}
\label{lemma:x1-solution}
Let $\mathfrak h$ be defined in \eqref{h1}. If $\tilde c\ne 0$, the update in~\eqref{eq:updatex_k} becomes 
\vspace{-0.1in}$$x_1^{k+1}=(\mathbf s_1 + \mathbf s_2)  T(\tilde c)/\|T(\tilde c)\|_2,
$$ where
$ \mathbf s_1 = \sqrt[3]{\frac{\mathbf c}{2 \mathfrak l_1} + \sqrt{\frac{1}{27}+\frac{\mathbf c^2}{4\mathfrak l_1^2}}}$, \,\, $\mathbf s_2 =\sqrt[3]{\frac{\mathbf c}{2\mathfrak l_1} - \sqrt{\frac{1}{27}+\frac{\mathbf c^2}{4\mathfrak l_1^2}}}$, \,\, $\mathbf c=\sqrt{\sum_{i=1}^d (|\tilde c_i|-\lambda_1)^2_+}$,  \,\,
 and \,\,
 \mbox{$T(\tilde c)=-(|\tilde c|-\lambda_1)_+ {\rm sgn}(\tilde c).$}
  If $\tilde c= 0$, the solution of \eqref{eq:updatex_k} is 
 $x_1^{k+1}=(\mathbf s_1 + \mathbf s_2) x'_1, $
 where $ \mathbf s_1$ and $ \mathbf s_2$ are determined as above with $\mathbf c = -\lambda_1$, and $x_1'$ is any vector with one component being 1 and the remaining components being 0.  
\end{lemma}
\proof{Proof.}
    See Appendix~\ref{lemma:solution}.
\QED

\textit{Update of block $x_2$.} Fix $x_1$, $y$ and $\omega$. Note that $x_2\mapsto \nabla_{x_2} \mathfrak f_i(x)=\omega_i \mathbf a_i + \beta(\phi_i(x)-y_i)\mathbf  a_i $ is Lipschitz continuous with constant $\beta\|\mathbf a_i\|^2 $. Hence, $x_2\mapsto \varphi_\beta(x,y,\omega)$ is $\mathfrak  l_2$-smooth, with $\mathfrak l_2=\beta\sum_{i=1}^q \|\mathbf a_i\|^2$. Thus we use the Lipschitz gradient surrogate for updating $x_2$:
\vspace{-0.1in}$$
\begin{array}{ll}
x_2^{k+1} &= \argmin_{x_2} \lambda_2 \|x_2\|_1 + \langle \nabla_{x_2}\varphi_\beta(x_1^{k+1},x_2^k,x_3^k,y^k,\omega^k),x_2-x_2^k\rangle + \frac{\mathfrak l_2}{2}  \|x_2-x_2^k\|_2^2\\
 &= \prox_{\frac{\lambda_2}{\mathfrak l_2}\|\cdot\|_1}\Big( x_2^k- \frac{1}{\mathfrak l_2} \nabla_{x_2}\varphi_\beta(x_1^{k+1},x_2^k,x_3^k,y^k,\omega^k) \Big). 
\end{array} $$

\textit{Update of block $x_3$.} Similarly to the update of $x_2$, we use the Lipschitz gradient surrogate:
$
x_3^{k+1} =  x_3^k-\frac{1}{\mathfrak l_3}  \nabla_{x_3}\varphi_\beta(x_1^{k+1},x_2^{k+1},x_3^k,y^k,\omega^k),
$
where $\mathfrak l_3=\beta q$. 

As noted in Remark~\ref{remark1}, Condition \ref{condition-C3} and Condition \ref{condition-C1} are satisfied.

\textit{Update of $y$.} We have 
\vspace{-0.1in}
$$\nabla h(y)=-\frac{1}{q}\Big[\frac{\mathbf b_i e^{-\mathbf b_i y_i}}{1+e^{-\mathbf b_i y_i}}\Big|_{i\in[q]} \Big], \quad \nabla^2 h(y)=\frac{1}{q}{\rm Diag}\Big(\frac{\mathbf b_i^2 e^{-\mathbf b_i y_i} }{(1+e^{-\mathbf b_i y_i})^2} \Big|_{i\in[q]}\Big)  \preceq \frac{1}{4q} I. 
$$
Hence $L_h= \frac{1}{4q}$.  The update of $y$ is as in \eqref{eq:y_update} (note that $\mathcal B=-I$ and $b=0$).

\textit{Choose $\beta$.}
We choose $\beta\geq 10L_h$ to satisfy \eqref{para_condition_2}.

\subsection{Applying the method in \cite{Drusvyatskiy2019}.}
Problem \eqref{svm_quad} can also be solved by the prox-linear method proposed in \cite{Drusvyatskiy2019}  
 \begin{equation}
\label{eq:prox-linear}
x^{k+1}\in \argmin_x  g(x) + \tilde h(\tilde \phi(x^k) + \nabla \tilde \phi(x^k)^\top(x-x^k))+ \frac{1}{2\tau}\|x-x^k\|^2,
\end{equation} 
where $\tilde h(y)=1/q\sum_{i=1}^q\log(1+e^{-y_i})$, $\tilde \phi_i(x_1,x_2,x_3)=\mathbf b_i(\langle \mathbf a_i,x_1 \rangle^2 + \langle \mathbf a_i, x_2\rangle  + x_3)$ and $\tau=(L_{\tilde h} L_{\tilde \phi})^{-1}$. 
 The subproblem~\eqref{eq:prox-linear} does not have a closed-form solution, we thus use the accelerated proximal gradient method \cite{Nesterov2004} to solve it approximately.

 \subsection{Numerical results for synthetic data sets.} 
\begin{figure}
    \centering
    \includegraphics[scale=0.3]{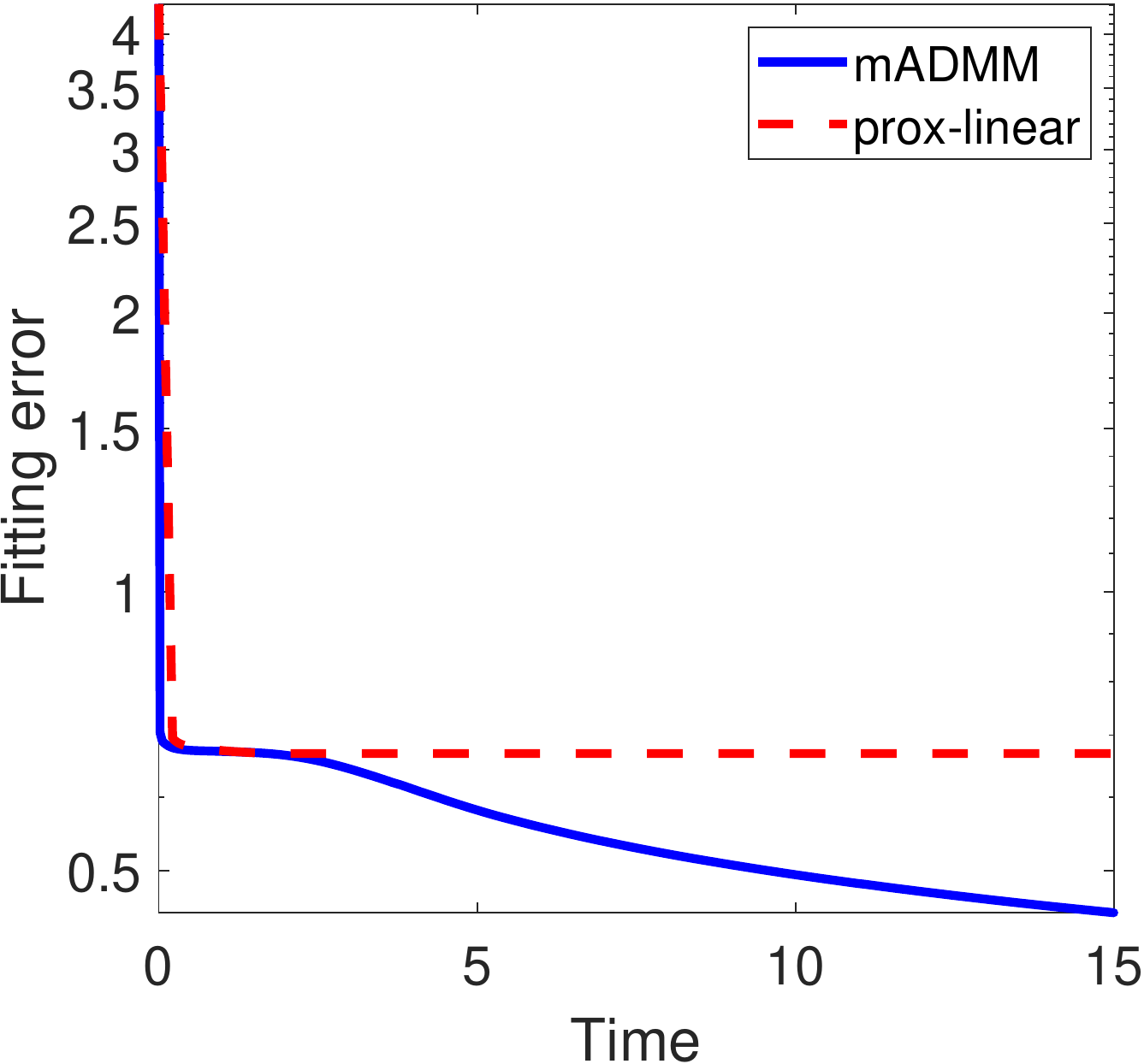} 
     \includegraphics[scale=0.3]{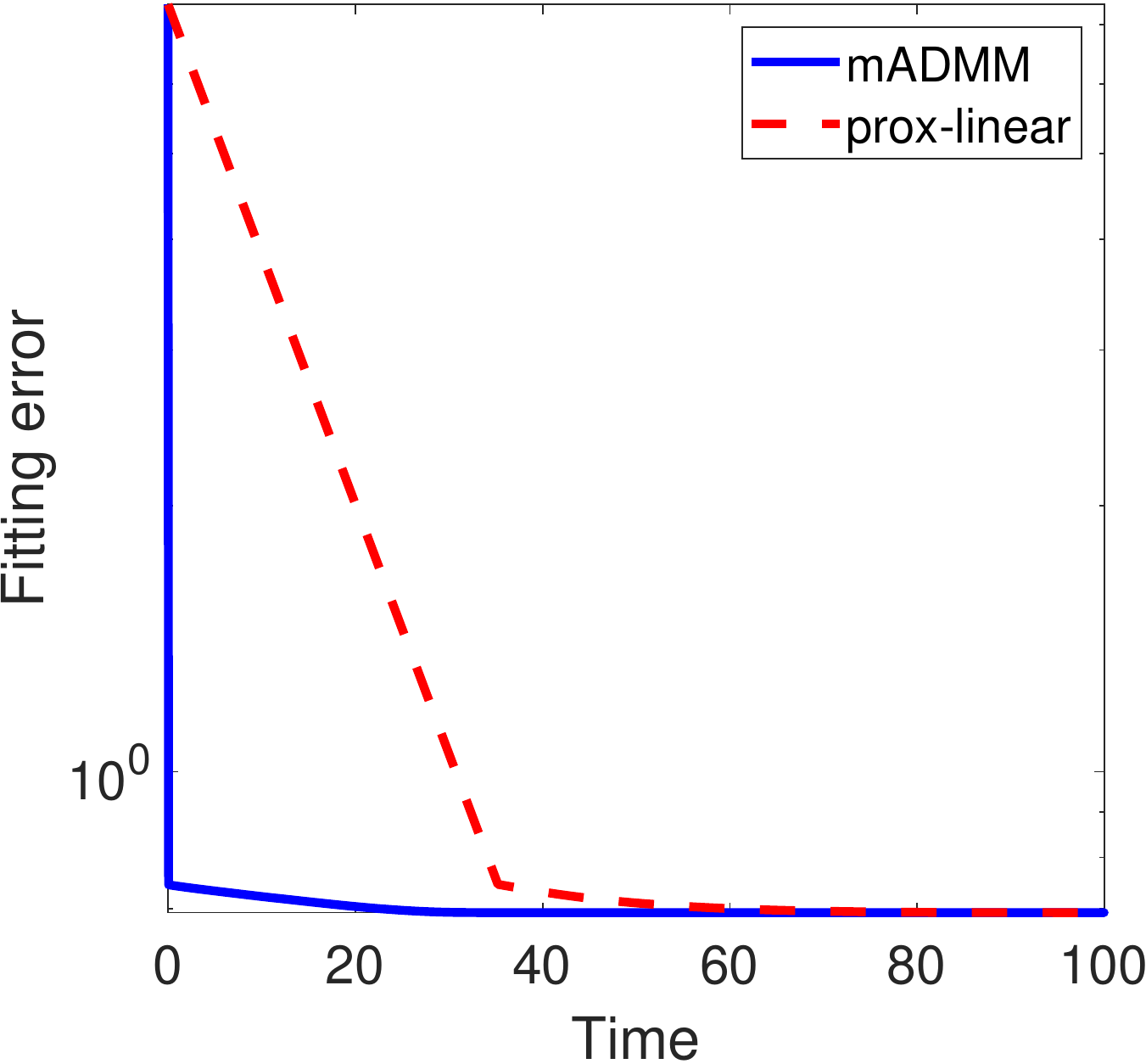} 
      \includegraphics[scale=0.312]{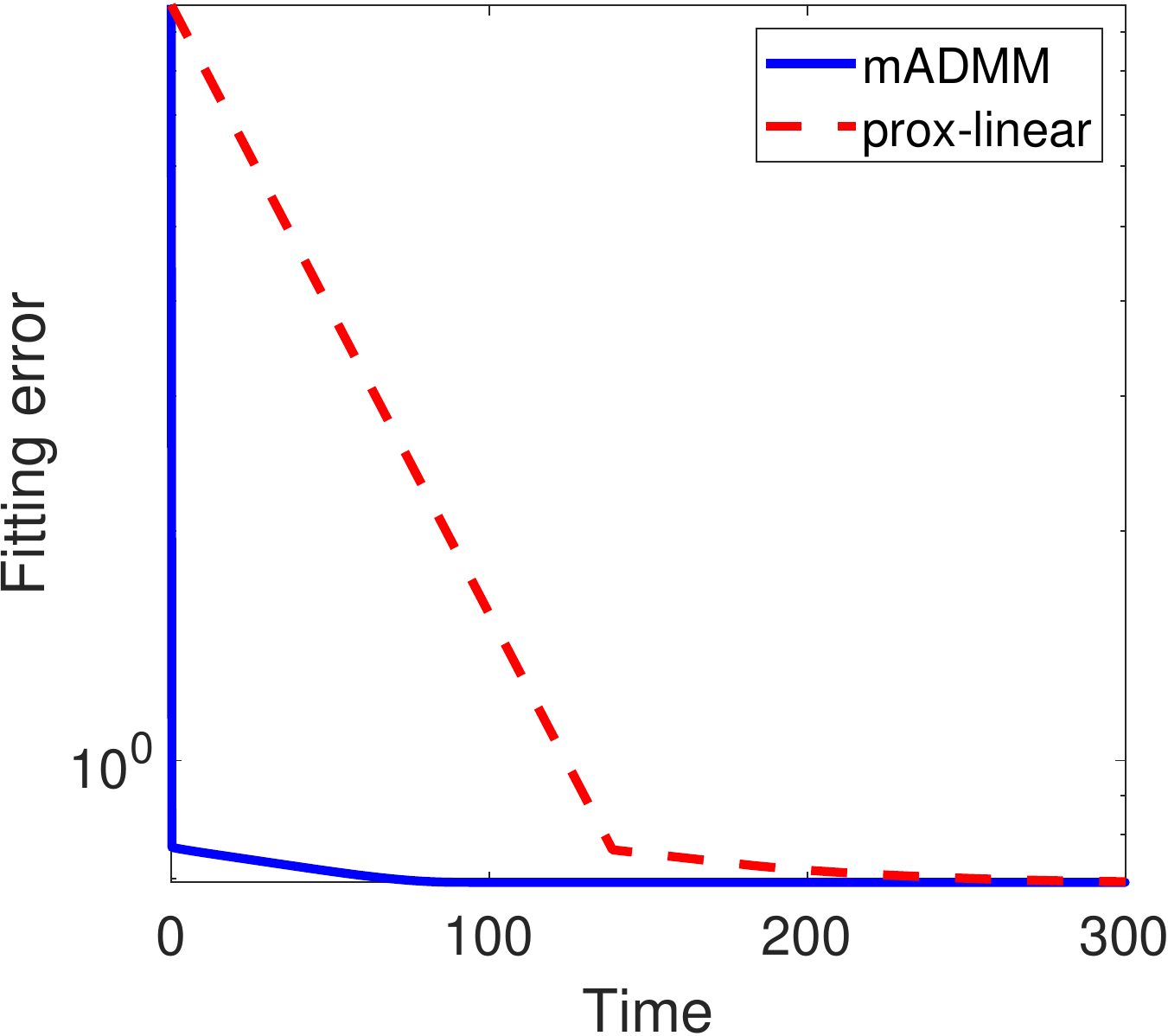} 
    \caption{Evolution of the fitting error with respect to the running time for 3 synthetic data sets with $(d,q)= (1\,000,100)$ (left), $(d,q)=(5\,000,1\,000)$ (middle) and $(d,q)=(10\,000,2\,000)$ (right).
     \label{fig:synthetic}} 
    \end{figure}

  We take $\lambda_1=0.001$ and $\lambda_2=0.1$.
 We use the Matlab command \texttt{rand} to generate a random matrix $A\in \bbR^{d\times q}$ (the columns of $A$ are $\mathbf a_i$, $i\in[q]$) and   \texttt{randsample} to generate a random vector $\mathbf b\in \bbR^q$, $\mathbf b_i\in \{-1,1\}$. The columns of $A$ are then normalized.  For \mADMM, we take  $\beta=2.5/q$. We run the algorithms with the same initial point generated by \texttt{rand} and with the same running time 15 seconds, 100 seconds and 300 seconds for the size $(d,q)= (1000,100), (5000,1000)$ and $(10\,000, 5000)$, respectively. We have tried different random initial points and observe that the results are similar. We report the evolution of the fitting error (i.e., the objective function of \eqref{svm_quad}) with respect to time in Figure~\ref{fig:synthetic} and the final fitting error in Table~\ref{tab:logistic-synthetic}.
        \begin{table}[t]
\centering
\caption{The final fitting error in solving Problem~\eqref{svm_quad} with synthetic data sets. \label{tab:logistic-synthetic}}
\begin{tabular}{|c|c|c|}
   \hline
$(d,q)$ & \mADMM & prox-linear \\
\hline
$(1\,000,100)$ & 0.450111 & 0.668748 \\
\hline
$(5000,1\,000)$ & 0.693019 & 0.693042\\
$(10\,000, 5\,000)$ & 0.692809 & 0.694131\\
 \hline 
  \end{tabular}
\end{table}
We observe from Figure \ref{fig:synthetic} and Table \ref{tab:logistic-synthetic} that \mADMM\, converges faster than prox-linear and \mADMM\, obtains better final fitting error than prox-linear.

 \subsection{Numerical results for real data sets.} 
\begin{figure}
    \centering
    \includegraphics[scale=0.3]{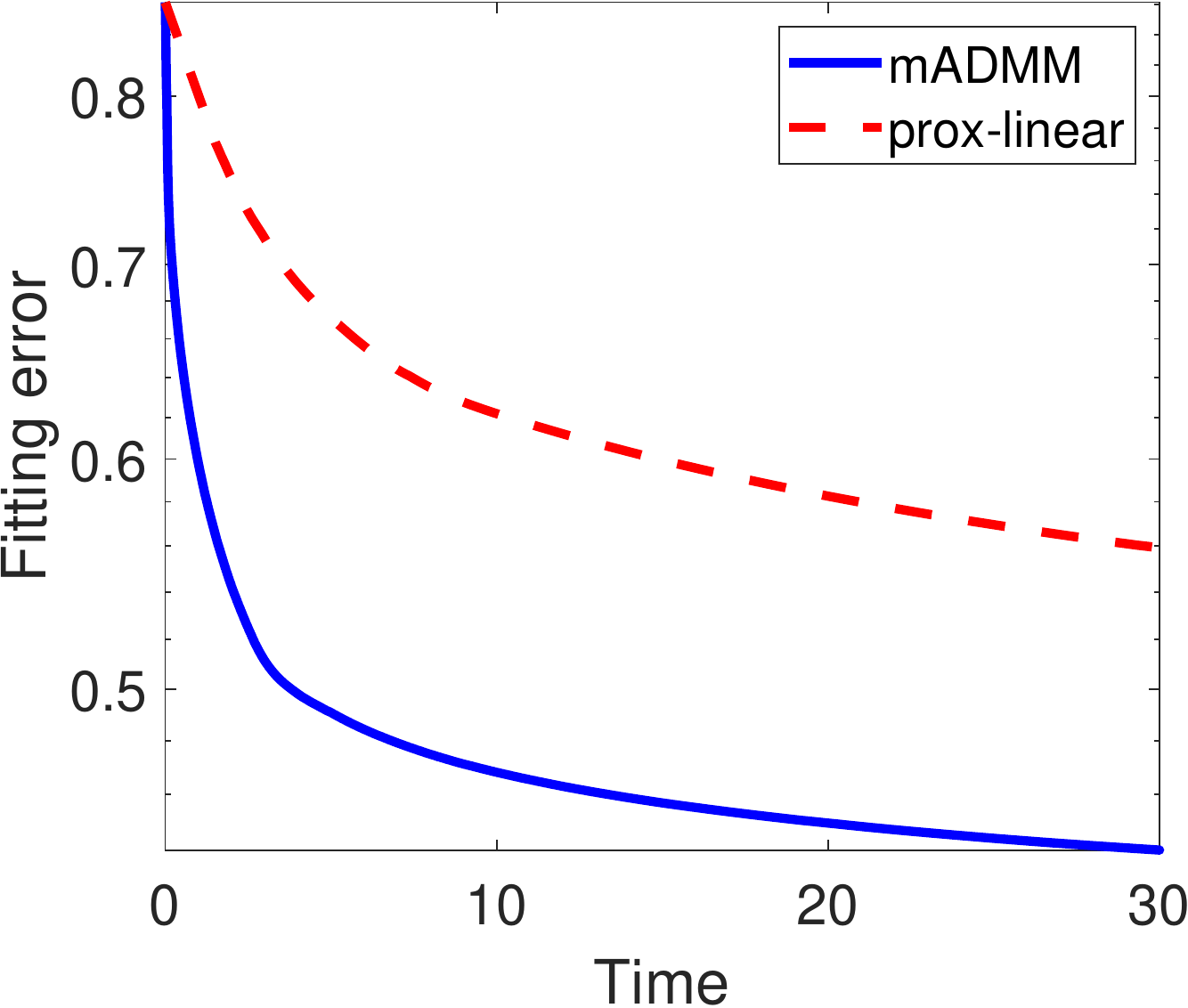} 
     \includegraphics[scale=0.3]{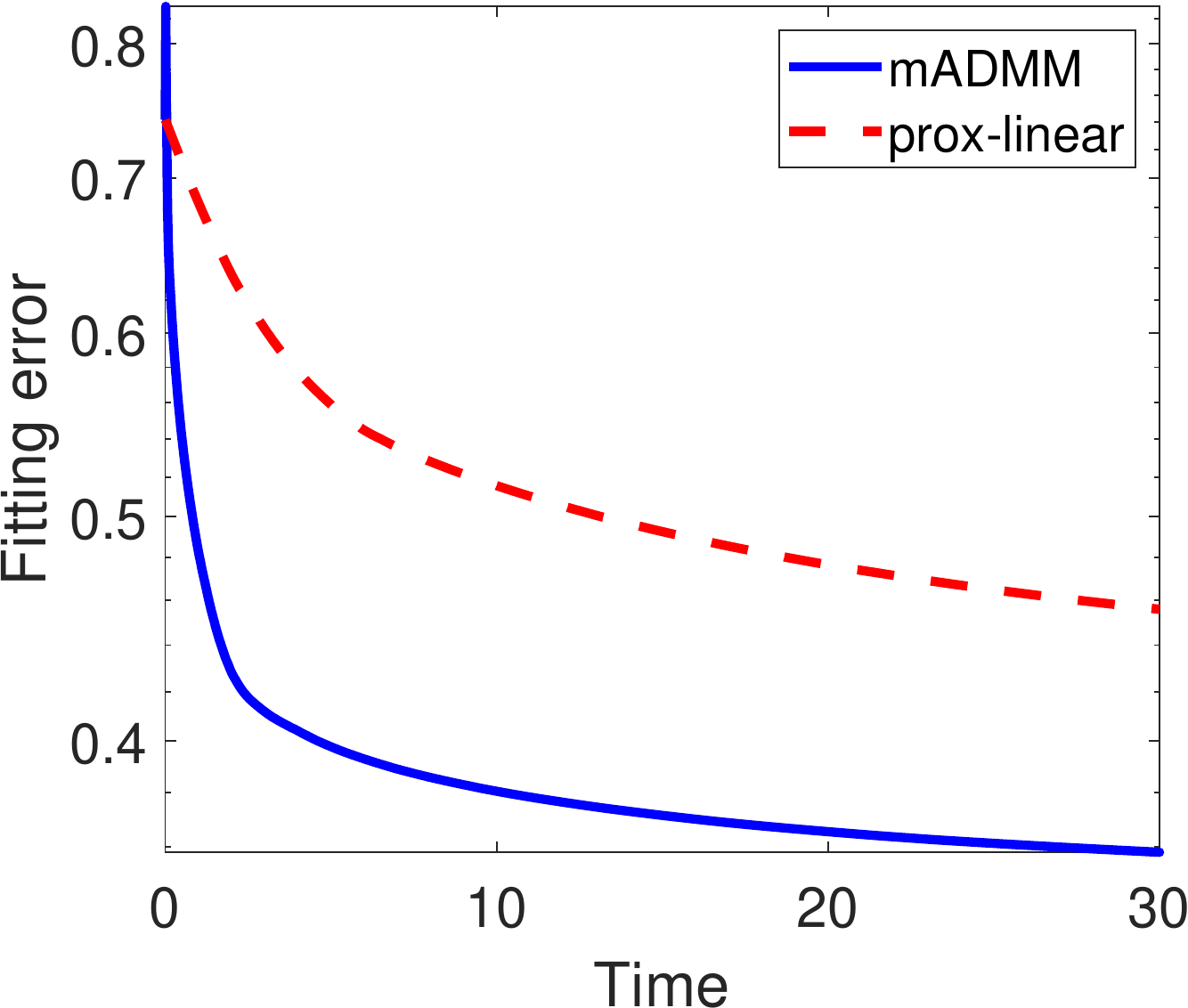} 
      \includegraphics[scale=0.3]{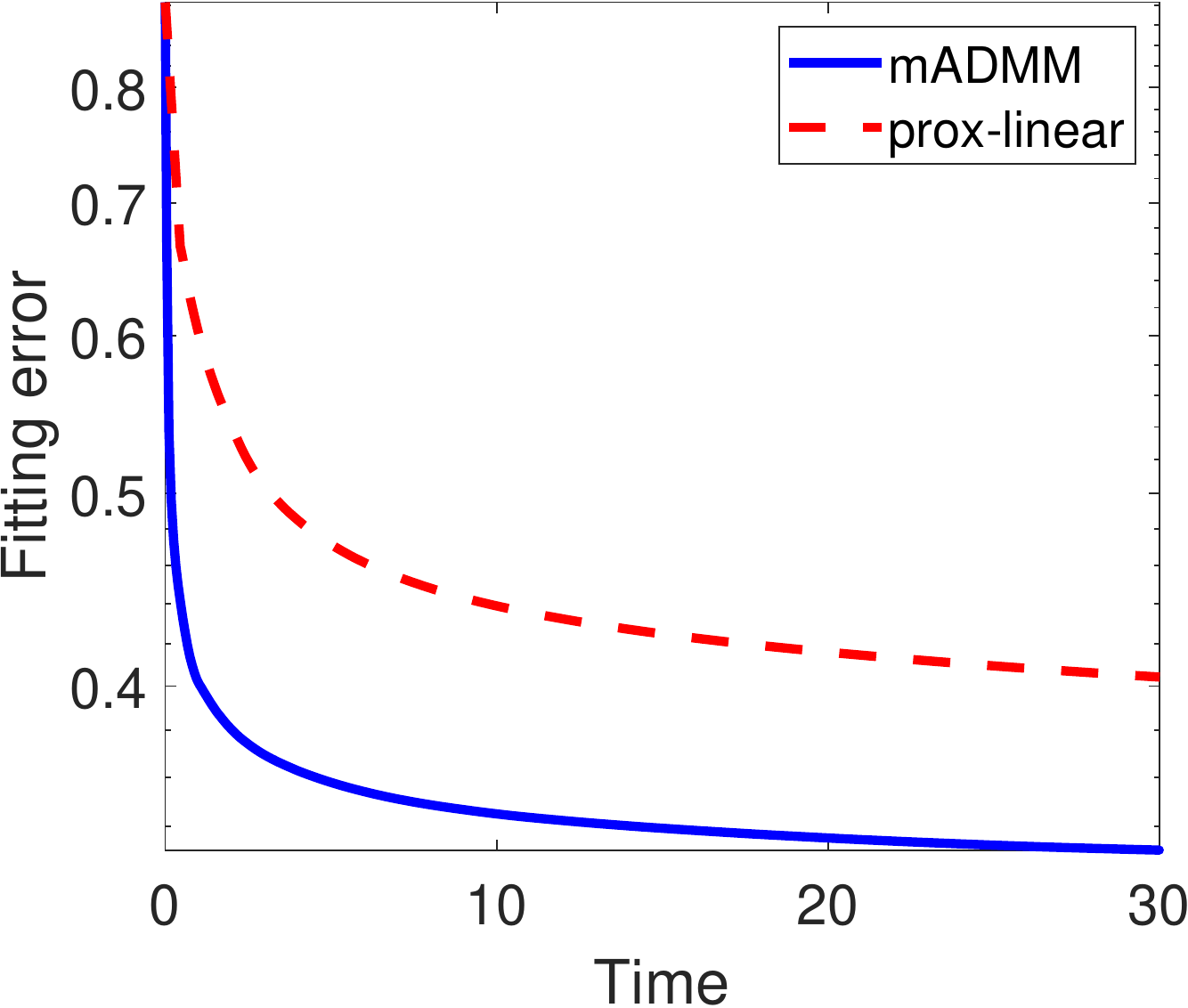} 
    \caption{Evolution of the fitting error with respect to the running time for 3 data sets duke breast-cancer (left), leukemia (middle) and colon-cancer (right).
     \label{fig:real}} 
    \end{figure}

In this experiment, we test the algorithms on three data sets leukemia, duke breast-cancer and colon-cancer. The data sets are available at \url{https://www.csie.ntu.edu.tw/~cjlin/libsvmtools/datasets/}. 
For each data set, we run each algorithm 30 seconds with the same random initial point  (note that we have tried many random initial points and the results are similar). We set $\lambda_1=\lambda_2=0.001$ and normalize the column of the instance matrix $A$ before running the algorithms. We report the evolution of the fitting error with respect to time in Figure~\ref{fig:real} and the final fitting error in Table~\ref{tab:real}.
    The results observed from Table~\ref{tab:real} and Figure~\ref{fig:real} are consistent with the results for synthetic data sets: \mADMM\, outperforms prox-linear, and, especially, \mADMM\, obtained better final fitting errors than prox-linear. 
  \begin{table}[t]
\centering
\caption{The final fitting error of \mADMM\, and prox-linear in solving~\eqref{svm_quad} with real data sets. 
\label{tab:real}}
\begin{tabular}{|c|c|c|c|}
   \hline
Data set&$(d,q)$ & \mADMM & prox-linear \\
\hline
duke breast-cancer &$(7\,129,44)$ & 0.440088 & 0.559264 \\
\hline
leukemia&$(7\,129,38)$ & 0.358154  & 0.455985\\
colon-cancer &$(2\,000,62)$ & 0.33082  & 0.404217\\
 \hline
  \end{tabular}
\end{table}

 \section{Conclusion}
 \label{sec:conclusion}
\vspace{-0.05in}
We have proposed \mADMM, a multiblock alternating direction method of multipliers, for solving a class of multiblock nonconvex optimization problems with nonlinear coupling constraints. Subsequential convergence, iteration complexity and global convergence are studied for the proposed method. One flexible feature of mADMM is that it allows the usage of block surrogate functions in updating the primal variables. The flexibility and the advantage of this feature are illustrated through the application of mADMM to solve an $l_1$-regularized logistic regression problem with a nonlinear classifier. By choosing suitable block surrogate functions, subproblems for updating the block variables have closed-form solutions. 
The numerical results have shown certain efficacy of \mADMM.    
We make an ending remark by proposing a potential future research direction, that is developing randomized/stochastic version of mADMM for solving Problem \eqref{nlq_2}.

\bibliographystyle{amsplain}
\bibliography{references}

\providecommand{\bysame}{\leavevmode\hbox to3em{\hrulefill}\thinspace}
\providecommand{\MR}{\relax\ifhmode\unskip\space\fi MR }
\providecommand{\MRhref}[2]{%
  \href{http://www.ams.org/mathscinet-getitem?mr=#1}{#2}
}
\providecommand{\href}[2]{#2}
\begin{thebibliography}{10}

\bibitem{Attouch2010}
H.~Attouch, J.~Bolte, P.~Redont, and A.~Soubeyran, \emph{Proximal alternating
  minimization and projection methods for nonconvex problems: An approach based
  on the {K}urdyka-{{\L}}ojasiewicz inequality}, Mathematics of Operations
  Research \textbf{35} (2010), no.~2, 438--457.

\bibitem{Attouch2013}
H.~Attouch, J.~Bolte, and B.~F. Svaiter, \emph{Convergence of descent methods
  for semi-algebraic and tame problems: proximal algorithms, forward--backward
  splitting, and regularized gauss--seidel methods}, Mathematical Programming
  \textbf{137} (2013), no.~1, 91--129.

\bibitem{Attouch2009}
Hedy Attouch and J{\'e}r{\^o}me Bolte, \emph{On the convergence of the proximal
  algorithm for nonsmooth functions involving analytic features}, Mathematical
  Programming \textbf{116} (2009), no.~1, 5--16.

\bibitem{Bolte2014}
J.~Bolte, S.~Sabach, and M.~Teboulle, \emph{Proximal alternating linearized
  minimization for nonconvex and nonsmooth problems}, Mathematical Programming
  \textbf{146} (2014), no.~1, 459--494.

\bibitem{Bolte2018_ADMM}
Jérôme Bolte, Shoham Sabach, and Marc Teboulle, \emph{Nonconvex
  {L}agrangian-based optimization: Monitoring schemes and global convergence},
  Mathematics of Operations Research \textbf{43} (2018), no.~4, 1210--1232.

\bibitem{Emilie2016}
E.~Chouzenoux, J.-C. Pesquet, and A.~Repetti, \emph{A block coordinate variable
  metric forward–backward algorithm}, Journal of Global Optimization
  \textbf{66} (2016), 457--485.

\bibitem{Cohen2021}
Eyal Cohen, Nadav Hallak, and Marc Teboulle, \emph{A dynamic alternating
  direction of multipliers for nonconvex minimization with nonlinear functional
  equality constraints}, Journal of Optimization Theory and Applications
  \textbf{193} (2022), 324–353.

\bibitem{Drusvyatskiy2019}
D.~Drusvyatskiy and C.~Paquette, \emph{Efficiency of minimizing compositions of
  convex functions and smooth maps}, Mathematical Programming \textbf{178}
  (2019), 1436--4646.

\bibitem{Duchi2018}
John~C Duchi and Feng Ruan, \emph{{Solving (most) of a set of quadratic
  equalities: composite optimization for robust phase retrieval}}, Information
  and Inference: A Journal of the IMA \textbf{8} (2018), no.~3, 471--529.

\bibitem{Dutter1981}
Rudolf Dutter and Peter~J. Huber, \emph{Numerical methods for the nonlinear
  robust regression problem}, Journal of Statistical Computation and Simulation
  \textbf{13} (1981), no.~2, 79--113.

\bibitem{HienNicolas_KLNMF}
Le~Thi~Khanh Hien and Nicolas Gillis, \emph{Algorithms for nonnegative matrix
  factorization with the {K}ullback-{L}eibler divergence}, Journal of
  Scientific Computing \textbf{87} (2021).

\bibitem{Hien2021_iADMM}
Le~Thi~Khanh Hien, Duy~Nhat Phan, and Nicolas Gillis, \emph{Inertial
  alternating direction method of multipliers for non-convex non-smooth
  optimization}, Computational Optimization and Applications \textbf{83}
  (2022), 247--285.

\bibitem{Titan2020}
\bysame, \emph{An inertial block majorization minimization framework for
  nonsmooth nonconvex optimization}, 2023.

\bibitem{Hien2022}
Le~Thi~Khanh Hien, Duy~Nhat Phan, Nicolas Gillis, Masoud Ahookhosh, and
  Panagiotis Patrinos, \emph{Block {Bregman} majorization minimization with
  extrapolation}, SIAM Journal on Mathematics of Data Science \textbf{4}
  (2022), no.~1, 1--25.

\bibitem{Kecman2005}
V.~Kecman, \emph{Support vector machines -- an introduction}, pp.~1--47,
  Springer Berlin Heidelberg, Berlin, Heidelberg, 2005.

\bibitem{Lu2015}
C.~{Lu}, J.~{Tang}, S.~{Yan}, and Z.~{Lin}, \emph{Nonconvex nonsmooth low rank
  minimization via iteratively reweighted nuclear norm}, IEEE Transactions on
  Image Processing \textbf{25} (2016), no.~2, 829--839.

\bibitem{lu2018relatively}
Haihao Lu, Robert~M. Freund, and Yurii Nesterov, \emph{Relatively smooth convex
  optimization by first-order methods, and applications}, SIAM Journal on
  Optimization \textbf{28} (2018), no.~1, 333--354.

\bibitem{Luss2013}
Ronny Luss and Marc Teboulle, \emph{Conditional gradient algorithms for
  rank-one matrix approximations with a sparsity constraint}, SIAM Review
  \textbf{55} (2013), no.~1, 65--98.

\bibitem{Maillard2010}
Sebastien Maillard, Thierry Roncalli, and Jerome Teiletche, \emph{The
  properties of equally weighted risk contribution portfolios}, The Journal of
  Portfolio Management \textbf{36} (2010), no.~4, 60--70.

\bibitem{melo2017}
Jefferson~G. Melo and Renato D.~C. Monteiro, \emph{Iteration-complexity of a
  {J}acobi-type non-{E}uclidean {ADMM} for multi-block linearly constrained
  nonconvex programs}, 2017.

\bibitem{Nesterov2004}
Y.~Nesterov, \emph{Introductory lectures on convex optimization: A basic
  course}, Springer New York, New York, 2004.

\bibitem{Ochs2019}
P.~Ochs, \emph{Unifying abstract inexact convergence theorems and block
  coordinate variable metric ipiano}, SIAM Journal on Optimization \textbf{29}
  (2019), no.~1, 541--570.

\bibitem{Boyd2014}
N.~Parikh and S.~Boyd, \emph{Proximal algorithms}, Foundations and Trends in
  Optimization \textbf{1} (2014), no.~3, 127--239.

\bibitem{Razaviyayn2013}
M.~Razaviyayn, M.~Hong, and Z.~Luo, \emph{A unified convergence analysis of
  block successive minimization methods for nonsmooth optimization}, SIAM
  Journal on Optimization \textbf{23} (2013), no.~2, 1126--1153.

\bibitem{Rock1981}
R.~Tyrrell Rockafellar, \emph{The theory of subgradients and its applications
  to problems of optimization - convex and nonconvex functions}, Heldermann,
  Heidelberg, Berlin, 1981.

\bibitem{RockWets98}
R.~Tyrrell Rockafellar and Roger J.-B. Wets, \emph{Variational analysis},
  Springer Verlag, Heidelberg, Berlin, New York, 1998.

\bibitem{Roosta2014}
Farbod Roosta-Khorasani, Kees van~den Doel, and Uri Ascher, \emph{Stochastic
  algorithms for inverse problems involving pdes and many measurements}, SIAM
  Journal on Scientific Computing \textbf{36} (2014), no.~5, S3--S22.

\bibitem{Wang2019}
Yu~Wang, Wotao Yin, and Jinshan Zeng, \emph{Global convergence of {ADMM} in
  nonconvex nonsmooth optimization}, Journal of Scientific Computing
  \textbf{78} (2019), 29--63.

\bibitem{Xu2013}
Y.~Xu and W.~Yin, \emph{A block coordinate descent method for regularized
  multiconvex optimization with applications to nonnegative tensor
  factorization and completion}, SIAM Journal on Imaging Sciences \textbf{6}
  (2013), no.~3, 1758--1789.

\bibitem{Xu2017}
\bysame, \emph{A globally convergent algorithm for nonconvex optimization based
  on block coordinate update}, Journal of Scientific Computing \textbf{72}
  (2017), no.~2, 700--734.

\end{thebibliography}

\begin{APPENDIX}{}
\section{Technical proofs}
\label{appendix2}

\subsection{Proof of Theory~\ref{thm:subsequential}.}
\label{thrm-sub}
The proof is similar to \cite[Theorem 1]{Hien2021_iADMM}. It follows from Proposition~\ref{prop:delta_converge} that if $\{(x^{k_n},y^{k_n},\omega^{k_n})\}$ converges to $(x^*,y^*,\omega^*)$  then  $\{(x^{k_n+1},y^{k_n+1},\omega^{k_n+1})\}$ and  $\{(x^{k_n-1},y^{k_n-1},\omega^{k_n-1})\}$ also converge to $(x^*,y^*,\omega^*)$.  On the other hand, from \eqref{eq:x_iupdate-2}, we have
\begin{align}
\label{Sub4}
u_i(x_i^{k+1},x^{k,i-1},y^k,\omega^k)+ g_i(x_i^{k+1})   \leq  u_i(x_i,x^{k,i-1},y^k,\omega^k) + g_i(x_i),\, \, \forall\,  x_i.
\end{align}
Choose $x_i=x_i^*$ in \eqref{Sub4} and note that $ u_i(x_i,z)$ is continuous, we have 
$
\limsup_{n\to\infty}   g_i(x_i^{k_n}) \leq   g_i(x_i^*). 
$
Furthermore, $g_i(x_i)$ is l.s.c. Hence, $ g_i(x_i^{k_n})\to g_i(x_i^*)$. 
Let $k=k_n\to\infty$ in~\eqref{Sub4}, for all $x_i$ we have
\vspace{-0.1in} 
\begin{equation}
\label{eq:blah}
\begin{array}{ll}
\varphi_\beta(x^*,y^*,\omega^*) + g_i(x_i^*) &\leq  u_i(x_i,x^*,y^*,\omega^*) +g_i(x_i)\\
&\leq \varphi_\beta((x_i,x^*_{\ne i}),y^*,\omega^*) +  \bar e_i(x_i,x^*,y^*,\omega^*) + g_i(x_i). 
\end{array}
\end{equation}  
Hence, we have
$x_i^*\in\arg\min_{x_i}\varphi_\beta((x_i,x^*_{\ne i}),y^*,\omega^*) + \bar  e_i(x_i,x^*,y^*,\omega^*) + g_i(x_i)
$ since $\bar e_i(x^*_i,x^*,y^*,\omega^*)=0$.
Thus, 
$0 \in \partial_{x_i}\Big( \varphi_\beta(x^*,y^*,\omega^*) + \bar  e_i(x^*_i,x^*,y^*,\omega^*) + g_i(x^*_i)\Big)$. Furthermore, $\nabla_{x_i} \bar e_i(x^*_i,x^*,y^*,\omega^*)=0$. Hence, we have $0 \in \partial_{x_i} \mL_\beta(x^*,y^*,\omega^*)$. 
Similarly, we can prove that $0 \in \partial_{y} \mL(x^*,y^*,\omega^*) $. Moreover, we have $\Delta \omega^k= \omega^{k} - \omega^{k-1}= \beta (\phi(x^k) + \mcalB y^k) \to 0.$ Hence, $\partial_\omega \mL(x^*,y^*,\omega^*) = \phi(x^*) + \mcalB y^*=0.
$
Finally, since we assume $\partial F(x)=\partial_{x_1} F(x) \times \ldots \times \partial_{x_m} F(x)$, we have 
\begin{equation*}
\begin{array}{ll}
\partial \mL_\beta(x,y,\omega) &= \partial F(x)+ \nabla \Big(h(y) +  \langle \omega,\phi(x) +\mcalB y-b \rangle + \frac{\beta}{2} \|\phi(x) + \mcalB y\|^2\Big)\\
&=\partial_{x_1} \mL_\beta(x,y,\omega) \times\ldots\times \partial_{x_m} \mL(x,y,\omega) \times \partial_{y} \mL_\beta(x,y,\omega)\times \partial_{\omega} \mL(x,y,\omega). 
\end{array}
\end{equation*}
We conclude that $ 0\in \partial \mL_\beta(x^*,y^*,\omega^*)$. \QED
\subsection{Proof of Proposition \ref{prop:bounded-sequence}.}
\label{bounded-sequence}
Since $\sigma_{\mB}>0$ we have $\mB$ is surjective. On the other hand, since $ran\, \phi(x) \subseteq Im(\mB)$ we have there exist $\bar y^k$ such that $\phi(x^k) +\mB\bar y^k=0 $. We consider $k\geq 1$. 
   Now we have
  \begin{equation}
 \label{temp4}
 \begin{array}{ll}
 \mL^k &=F(x^{k})+ h(y^{k}) +\frac{\beta}{2}\|\phi(x^{k})+\mB y^{k}\|^2 +\langle \omega^k, \phi(x^{k})+\mB y^{k}\rangle\\
 &=F(x^{k}) + h(y^{k}) +\frac{\beta}{2}\|\phi(x^{k})+\mB y^{k}\|^2 + \langle \mB^*\omega^k,y^{k}-\bar y^k\rangle.
 \end{array}
 \end{equation} 
On the other hand, from \eqref{eq:y_update} we have 
 \begin{equation}
\label{eq:omegabound} \nabla h(u^k) + \mB^* \big(\omega^{k} + \beta(\phi(x^{k+1}) + \mB y^{k+1})\big) + L_h(y^{k+1}-y^k) =0.
\end{equation} 
Hence, \vspace{-0.1in}
$$
\begin{array}{ll}
\langle \mB^*\omega^k,y^{k}-\bar y^k\rangle
&= \big\langle \nabla h( y^k) + L_h \Delta y^{k+1} +   \mB^*\Delta w^{k+1}, \bar y^k-y^{k}\big\rangle\\
&\geq 
\langle\nabla h(y^k) , \bar y^k-y^{k}\rangle-\big( L_h \| \Delta y^{k+1}\| + \|\mB^*\Delta\omega^{k+1}\|\big) \| \bar y^k-y^{k}\|.
\end{array} 
$$
Together with~\eqref{temp4} and $L_h$-smooth property of $h$ we imply that
\vspace{-0.1in} \begin{equation}
\label{temp5}
\mL^k \geq F(x^{k}) + h(\bar y^k) - \frac{L_h}{2}\|y^k-\bar y^k\|^2- \big( L_h \| \Delta y^{k+1}\|  + \|\mB^* \Delta\omega^{k+1}\|\big) \| \bar y^k-y^{k}\|. 
\end{equation} 
Moreover, we have 
 \begin{equation}
\label{temp6}
\begin{split}
\| \bar y^k-y^{k}\|^2 \leq \frac{1}{\lambda_{\min}(\mB^*\mB)} \|\mB(\bar y^k-y^{k})\|^2&= \frac{1}{\lambda_{\min}(\mB^*\mB)}\|\phi( x^k)  + \mcalB y^{k} \|^2
\\
&=\frac{1}{\lambda_{\min}(\mB^*\mB)} \big\|\frac{1}{\beta } \Delta \omega^k\big\|^2.
\end {split}
\end{equation}
On the other hand, Proposition~\ref{prop:delta_converge} shows that $\|\Delta \omega^k\|$,  $\|\Delta x^k\|$ and $\|\Delta y^k\|$ converge to 0, and from~\eqref{recursive-2} we have $\mL^k$ is upper bounded. Therefore,  \eqref{temp5} and \eqref{temp6} imply that $ F(x^{k}) + h(\bar y^k)$ is upper bounded. So $\{x^k\}$ is bounded. Consequently, $\phi(x^k) $ is bounded.  
On the other hand, we have
$$ \|y^k\|^2\leq \frac{1}{\lambda_{\min}(\mB^*\mB)} \|\mB y^k\|^2 = \frac{1}{\lambda_{\min}(\mB^*\mB)}\big \|\frac{1}{\beta } \Delta \omega^k-\phi(x^k) \big\|^2.
$$
Hence, $\{y^k\}$ is bounded, which implies $\|\nabla h( y^k)\|$ is bounded. Finally, from~\eqref{eq:omegabound}, we have  $\{\omega^k\}$ is bounded. \QED

\subsection{Proof of Theorem \ref{thrm:global}.}
\label{proof-global}
We  do the analysis in the bounded set containing the generated sequence of \mADMM. 
The Lyapunov sequence 
$\tilde L_\beta(z^{k},y^{k-1})
$ has the following properties.

\emph{(i) Sufficient decreasing property.}
We derive from \eqref{recursive-2} that 
$$\mL^{k+1} + \frac12\tilde\eta \|\Delta x^{k+1}\|^2+ \frac{3\tilde\delta\hat{\delta}^2}{\beta\sigma_\mcalB} \| \Delta y^{k+1}\|^2 
\leq \mL^k+   \frac{3  \hat\delta^2}{\beta\sigma_\mcalB}\| \Delta y^{k}\|^2,
$$ 
where $\tilde\eta=\min_{i\in [m]} \underline{\eta}_i$. 
Hence
$$\tilde{\mL}_\beta (z^{k+1}, y^{k}) + \frac12\tilde\eta \|\Delta x^{k+1}\|^2+ 3\frac{(\tilde\delta -1) \hat\delta^2}{\beta\sigma_\mcalB} \| \Delta y^{k+1}\|^2 \leq \tilde{\mL}_\beta (z^{k}, y^{k-1}).
$$

\emph{(ii) Boundedness of subgradient.} We have 
\vspace{-0.1in}
$$
\begin{array}{ll}
\partial_{x_i} \tilde{\mL}_\beta (z^{k+1},y^{k})
&=\partial_{x_i} F(x^{k+1}) + \nabla_{x_i}\phi(x^{k+1})\big(\omega^{k+1} +\beta (\phi(x^{k+1})+\mcalB y^{k+1})\big)\\
\nabla_{\omega} \tilde{\mL}_\beta (z^{k+1}, y^{k})&=\phi(x^{k+1})+\mcalB y^{k+1},\\
 \nabla_y\tilde{\mL}_\beta (z^{k+1},y^{k})&= \nabla h(y^{k+1})+\mcalB^*\big(\omega^{k+1} +\beta (\phi(x^{k+1})+\mcalB y^{k+1})\big)+   \frac{6\hat\delta^2}{\beta\sigma_\mcalB}  (y^{k+1}- y^{k}),\\
  \nabla_{\tilde y}\tilde{\mL}_\beta (z^{k+1},y^{k}) &= \frac{6\hat\delta^2}{\beta\sigma_\mcalB} (y^{k}-y^{k+1}).
\end{array}  
$$
On the other hand, Proposition~\ref{prop:iteration_complexity} showed that there exists $\chi_i^k\in\partial_{x_i} F(x^k)$ such that \eqref{R_i_rate} holds. Therefore, it is not difficult to prove that there exist $\tilde\chi^{k+1}\in\partial \tilde{\mL}_\beta (x^{k+1},y^{k+1},\omega^{k+1},y^{k})  $ such that 
$\|\tilde\chi^{k+1}\|\leq a_1 \|\Delta x^{k+1}  \| + a_2 \|\Delta y^{k+1}  \|  + a_3 \|\Delta y^{k}  \| 
$ 
for some positive constants $a_1$, $a_2$ and $a_3$. 

\emph{(iii) K{\L}~property.} We assume that $\tilde \mL$ has the K{\L}~property with constant $\sigma_{\tilde \mL}$. 

\emph{(iv) A continuity property.} Suppose a subsequence $(z^{k_n}) \to(x^*,y^*,\omega^*).$ Proposition~\ref{prop:delta_converge} showed that $y^{k_n-1}\to y^*$.Moreover, in the proof of Theorem~\ref{thm:subsequential} we proved that $g_i(x^{k_n}_i)\to g_i(x_i^*).$ Therefore, $\mL_\beta(x^{k_n},y^{k_n},\omega^{k_n})\to \mL_\beta(x^*,y^*,\omega^*)$. Consequently, $\tilde{\mL}_\beta(z^{k_n},y^{k_n-1})\to \tilde{\mL}_\beta(x^*,y^*,\omega^*,y^*)$. 

We can prove 
$\sum\limits_{k=1}^\infty   \|\Delta x^{k+1}  \| + \|\Delta y^{k+1}  \|  + \|\Delta y^{k}  \| <\infty
$
by using the above properties and the same techniques of \cite[Theorem 1]{Bolte2014} (as this is typical technique, see e.g., \cite{Hien2021_iADMM,Xu2017}, we omit the details). Then $x^k\to x^*$ and $y^k\to y^*$.  
Moreover, from \eqref{R_i_rate}, we have $R^{k}_c = \frac{1}{\beta}\|\Delta \omega^k\|=O(\|\Delta y^k\| + \|\Delta y^{k-1}\|) $. Hence $\sum_{k=1}^\infty \|\Delta \omega^{k}  \|<\infty$, leading to $\omega^k \to\omega^*$. 
Finally, we use the same techniques of \cite[Theorem 2]{Attouch2009} 
 to obtain the convergence rate for $\{z^k\}_{k\geq 1}$. \QED

\subsection{Proof of Lemma \ref{lemma:x1-solution}.}  
\label{lemma:solution}
Suppose $\tilde c \ne 0$. We have \vspace{-0.1in}
  \begin{equation}
\label{num:1}
\begin{split}
\min_{x_1} \big\{\lambda_1 \|x_1\|_1 + \langle \tilde  c ,x_1\rangle : \|x_1\|_2= 1 \big\}&\eqa\min_{x_1} \big\{ \lambda_1\|x_1\|_1 + \langle \tilde  c ,x_1\rangle : \|x_1\|_2 \leq 1 \big\}\\
&\eqb -\sqrt{\sum_{i=1}^d (|\tilde c_i|-\lambda_1)^2_+}
\end{split}
\end{equation} 
and is solved by 
$x_1^*=T(\tilde c)/\|T(\tilde c)\|_2$. Here we used \cite[Proposition 4.6]{Luss2013} for  (b) and the fact $\|x_1^*\|_2=1$ for (a). 
 On the other hand, we have 
\begin{equation}
\label{nm1}
\begin{array}{ll}
&\min_{x_1\in \bbR^d}\lambda_1 \|x_1\|_1 + \langle \tilde  c ,x_1\rangle  + \frac{\mathfrak l_1}{4} \|x_1\|_2^4 + \frac{\mathfrak l_1}{2} \|x_1\|_2^2 \\
\equiv&\min_{x_1\in\bbR^d,t\in\bbR_+}\big\{  \lambda_1\|x_1\|_1 +\langle \tilde  c ,x_1\rangle + \frac{\mathfrak l_1}{4} t^4 +   \frac{\mathfrak l_1}{2} t^2 :  t^2= \|x_1\|_2^2  \big\}\\
\equiv&\min_{t\in\bbR_+} \Big \{\frac{\mathfrak l_1}{4} t^4 +   \frac{\mathfrak l_1}{2} t^2  + \min_{x_1}  \big\{ \lambda_1 \|x_1\|_1 +\langle \tilde  c ,x_1\rangle  :  \|x_1\|_2^2 = t^2  \big\} \Big\}\\
\equiv&  \min_{t\in\bbR_+}\Big \{  \frac{\mathfrak l_1}{4} t^4 +    \frac{\mathfrak l_1}{2} t^2 + t\min_{x_1'}   \big\{\lambda_1  \|x_1'\|_1 +\langle \tilde  c , x_1'\rangle  :  \|x_1'\|_2^2 = 1  \big\} \Big\}\\
\eqc& \min_{t\in\bbR_+}\Big \{   \frac{\mathfrak l_1}{4} t^4 +    \frac{\mathfrak l_1}{2} t^2 - t\mathbf c \Big\},
\end{array}
\end{equation}
where we have used \eqref{num:1} for (c). Note that $t^*=\mathbf s_1 + \mathbf s_2$ is the solution of the last minimization problem in \eqref{nm1}(which is the nonnegative real solution of the cubic equation $\mathfrak l_1(t^*)^3 + \mathfrak l_1 t^* -\mathbf c=0$). 

When $\tilde  c=0$, note that $\min_{x_1'}   \big\{\lambda_1  \|x_1'\|_1  :  \|x_1'\|_2^2 = 1  \big\}=\lambda_1$ and the optimal value is obtained at any point $x_1'$ that has only one component being 1 and the remaining components being 0. \QED

 \section{An example}
\label{appendix:example}
Suppose $f(x)=\sum_{i=1}^m \hat f_i(\|x_i\|_2), $ 
where $\hat f_i$ is a continuously differentiable concave function with Lipschitz gradient on any given bounded set. This covers many nonconvex regularizers of low rank representation problems, see e.g., \cite{Lu2015}.   Since $\hat f_i$ is concave, we have 
$$
\hat f_i(\| x_i\|_2) \leq \hat f_i(\|\tilde x_i\|_2) +  \nabla\hat f_i(\|\tilde x_i\|_2)( \|x_i\|_2 - \|\tilde x_i\|_2).
$$ 
Hence $u_i$ defined in the following is a block surrogate function of $\varphi_\beta$ 
\vspace{-0.1in} \begin{equation*}
u_i(x_i,\tilde x,y,\omega) = f(\tilde x_i,x_{\ne i}) +  \nabla\hat f_i(\|\tilde x_i\|_2)( \|x_i\|_2 - \|\tilde x_i\|_2)+ \hat u_i(x_i,\tilde x,y,\omega), 
\end{equation*}
where $\hat u_i(x_i,\tilde x,y,\omega)$ is a surrogate of $\hat\varphi_\beta(x,y,\omega) = \langle \omega , \phi(x) + \mcalB y\rangle + \frac{\beta}{2}\|\phi(x)+\mcalB y\|^2$ with respect to block $x_i$. Assume $\hat u_i$ is twice continuously differentiable and $\nabla_{x_i} \hat u_i(x_i,x,y,\omega)=\nabla_{x_i} \hat\varphi_\beta (x,y,\omega)$ for all $x,y,\omega$. Note that 
$$\partial_{x_i}  u_i(x_i,\tilde x,y,\omega)=  \nabla\hat f_i(\|\tilde x_i\|_2) \partial( \|x_i\|_2 ) + \nabla_{x_i} \hat u_i(x_i,\tilde x,y,\omega). $$
Hence, any subgradient in $\partial_{x_i}  u_i(x^{k+1}_i, x^{k,i-1},y^{k},\omega^{k})$ has the form 
$$S^k_i= \nabla\hat f_i(\| x^{k}_i\|_2)\xi^k_i +\nabla_{x_i} \hat u_i(x^{k+1}_i, x^{k,i-1},y^{k},\omega^{k}),
$$ 
where $\xi^k_i \in\partial (\|\cdot\|_2)(x_i^{k+1})$. 
Moreover, it follows from \cite[Corollary 5Q]{Rock1981} that 
$\partial \hat f_i(x_i) =\nabla f_i(\|x_i\|_2)\partial (\|x_i\|_2).
$
Thus we take 
$\bar S^k_i = \nabla f_i(\|x^{k+1}_i\|_2)\xi^k_i + \nabla_{x_i} \hat\varphi_\beta(x^{k+1},y^k,\omega^k) \in \partial_{x_i} \varphi(x^{k+1},y^{k},\omega^{k}).$
Assuming $(x^k,y^k,\omega^k)$ is bounded, we have
 \begin{align*}
&\| S^k_i- \bar S^k_i\|=\big\| \big(\nabla\hat f_i(\| x^k_i\|_2)- \nabla \hat f_i(\|x^{k+1}_i\|_2)\big) \xi^k_i + \nabla_{x_i} \hat u_i(x^{k+1}_i,x^{k,i-1},y^k,\omega^k) \\
&\qquad\qquad-\nabla_{x_i} \hat u_i(x^{k+1}_i,x^{k+1},y^k,\omega^k)\big\|\leq L_{\hat f_i} (\|x^{k+1}_i - x^k_i\|)\|\xi^k_i\| + \bar L'_{i} \|x^{k+1}-x^k\|,
\end{align*}
where $\bar L'_{i}$ and $L_{\hat f_i}$ are the Lipschitz constant of  $\nabla_{x_i} \hat u_i$ and $\nabla \hat f_i$ on the bounded set containing $(x^k,y^k,\omega^k)$. Then \eqref{eq:l1} is satisfied with $L_i=L_{\hat f_i} + \bar L_i'$.
 \end{APPENDIX}


\end{document}